\documentclass{amsart}

\usepackage{mathtools,amssymb,amsthm}
\usepackage[unicode,hidelinks]{hyperref}
\usepackage{enumitem}

\theoremstyle{definition}
\newtheorem{theorem}{Theorem}[section]
\newtheorem{lemma}[theorem]{Lemma}
\newtheorem{corollary}[theorem]{Corollary}
\newtheorem{conjecture}[theorem]{Conjecture}
\newtheorem{remark}[theorem]{Remark}
\newtheorem{example}[theorem]{Example}

\newcommand{\CC}{\mathbb{C}}

\newcommand{\QQ}{\mathbb{Q}}
\newcommand{\RR}{\mathbb{R}}
\newcommand{\ZZ}{\mathbb{Z}}
\newcommand{\QQalg}{\QQ^{\mathrm{alg}}}

\DeclareMathOperator{\reg}{reg}
\DeclareMathOperator{\Gal}{Gal}
\DeclareMathOperator{\Norm}{N}

\begin{document}

\title[Algebraic independence of regulators]{Algebraic independence of regulators of number fields Galois over the rational number field under the weak Schanuel conjecture}

\author[Y.~Iizuka]{Yoshichika Iizuka}
\address{Department of Mathematics, Gakushuin University, Mejiro, Toshima-ku, Tokyo 171-8588, Japan}
\email{iizuka@math.gakushuin.ac.jp}

\author[Y.~Konomi]{Yutaka Konomi}
\address{Department of Mathematics, Meijo University, Tempaku-ku, Nagoya 468-8502, Japan}
\email{konomi@meijo-u.ac.jp}



\begin{abstract}
Under the weak Schanuel conjecture, we show that, for two number fields that are Galois over \(\QQ\),
the following two conditions are equivalent: the unit groups of both fields have positive rank and the maximal real subfields of the two fields are distinct;
the regulators of the two fields are algebraically independent over \(\QQ\).
Furthermore, under the same conjecture, we show that, for finitely many number fields that are Galois over \(\QQ\), if the unit groups of all these fields have positive rank
and the unit group of the intersection of any two distinct fields has rank \(0\), then the regulators of these fields are algebraically independent over \(\QQ\).
Under the same conjecture, we also show that, for two number fields that are Galois over \(\QQ\), the maximal real subfields of the two fields are distinct if and only if the ratio of their regulators is a transcendental number.
\end{abstract}

\maketitle

\section{Introduction}\label{sec:introduction}

In this paper, a number field means a subfield of the field \(\CC\) of complex numbers that is a finite extension of \(\QQ\).
An algebraic number is a complex number that is a root of a nonzero polynomial with rational coefficients, and an algebraic integer is a complex number that is a root of a monic polynomial with integer coefficients. A transcendental number is a complex number that is not an algebraic number.
We write \(\QQalg\) for the field of all algebraic numbers in \(\CC\).

The problem of characterizing number fields by their invariants is one of the fundamental problems in algebraic number theory.
The regulator is a fundamental invariant of a number field, alongside the discriminant and the class number.
The analytic class number formula relates these invariants to the residue of the Dedekind zeta function at \(s = 1\).
The question of the extent to which the regulator distinguishes number fields is closely related to the problem of distinguishing number fields by their Dedekind zeta functions.
For a number field \(K\), we write \(\reg(K)\) for its regulator and \(\zeta_{K}(s)\) for its Dedekind zeta function.
We say that two number fields \(K_{1}, K_{2}\) are arithmetically equivalent if \(\zeta_{K_{1}}(s) = \zeta_{K_{2}}(s)\) holds identically.
Perlis characterized arithmetic equivalence by Gassmann's group-theoretic condition
and showed, in particular, that if two number fields are arithmetically equivalent and one of them is Galois over \(\QQ\), then the two fields are isomorphic~\cite{Perlis1977}.

On the other hand, de Smit and Perlis showed that two number fields with the same Dedekind zeta function do not necessarily have the same class number~\cite{SmitPerlis1994}.
Two arithmetically equivalent number fields have the same numbers of real and complex embeddings, the same discriminant, and the same number of roots of unity~\cite{Perlis1977}.
Let \(K_{1}, K_{2}\) be two arithmetically equivalent number fields, and write \(h_{1}, h_{2}\) for their respective class numbers.
The analytic class number formula then gives
\[
h_{1}\reg(K_{1}) = h_{2}\reg(K_{2}).
\]
Hence the above result of de Smit and Perlis also shows that the Dedekind zeta function does not determine the regulator in general.

In the direction of strengthening arithmetic equivalence, Prasad introduced refined Gassmann triples~\cite{Prasad2017}.
This notion strengthens the condition that the permutation modules with rational coefficients are isomorphic
to the condition that the permutation modules with integer coefficients are isomorphic.
Prasad showed that two number fields arising from a refined Gassmann triple not only have the same Dedekind zeta function but also have isomorphic class groups and isomorphic idele class groups.
Furthermore, Sutherland introduced local integral equivalence as a notion that strengthens arithmetic equivalence,
and showed that two locally integrally equivalent number fields are arithmetically equivalent, have isomorphic class groups, and have the same regulator~\cite[Proposition~3.4]{Sutherland2021}.

Gun, Murty, and Rath discuss the transcendence of the regulator of a number field whose unit group has positive rank~\cite[proof of Theorem~5.1]{GunMurtyRath2011}.
In that paper, they express the regulator as the value of a polynomial with rational coefficients at the logarithms of the absolute values of conjugates of units,
and they use this expression to deduce the transcendence of the regulator under the weak Schanuel conjecture.
Recently, Cruz studied relations among the regulators of several totally real number fields under the same conjecture.
More precisely, Cruz showed that, in a finite family of totally real number fields different from \(\QQ\),
if the Dedekind zeta functions are pairwise distinct,
then the regulators of these fields are linearly independent over \(\QQalg\)~\cite[Theorem~5.3]{Cruz2025}.
Cruz also showed that, under the weak Schanuel conjecture, two totally real number fields
have the same regulator if and only if they have the same Dedekind zeta function and the same class number~\cite[Corollary~5.4]{Cruz2025}.
This result implies the following: under the weak Schanuel conjecture, if totally real number fields \(K_{1}, K_{2}\) that are Galois over \(\QQ\) satisfy \(K_{i} \ne \QQ\) \((i = 1, 2)\)
and have the same regulator, then the two fields coincide as subfields of \(\CC\).
Indeed, in this case the above result of Cruz shows that the Dedekind zeta functions also coincide,
and hence the result of Perlis yields an isomorphism \(\varphi \colon K_{1} \longrightarrow K_{2}\).
If we regard this \(\varphi\) as an embedding \(K_{1} \hookrightarrow \CC\) of \(K_{1}\), then the normality of \(K_{1}/\QQ\) gives
\[
K_{2} = \varphi(K_{1}) = K_{1}.
\]

In contrast, for totally imaginary number fields that are Galois over \(\QQ\), there exist distinct fields with the same regulator.
Iizuka and Konomi constructed an infinite family of pairs of distinct totally imaginary biquadratic fields
and an infinite family of pairs of distinct totally imaginary cyclic quartic fields~\cite[Theorems~1.2 and~1.3]{IizukaKonomi2026}.
In each of these infinite families, the two fields in each pair have the same regulator.
Moreover, in the same paper, they gave examples of pairs of distinct quartic CM fields that are Galois over \(\QQ\)
such that the two fields have the same discriminant, the same regulator, and the same class number~\cite[Examples~7.1 and~7.2]{IizukaKonomi2026}.
These examples show that, even for two number fields that are Galois over \(\QQ\) and whose unit groups have positive rank,
the coincidence of the regulators alone does not imply that the two fields coincide.

For a number field \(K\), we write \(r(K)\) for the rank of its unit group and \(K^{+}\) for its maximal real subfield \(K \cap \RR\).
We give the precise definitions of \(r(K)\) and of the regulator \(\reg(K)\) in Section~\ref{sec:preliminaries}.

The two main theorems below assume a special case of the weak Schanuel conjecture, whose precise formulation we give as Conjecture~\ref{conj:wsc} \textup{(WSC)} in Section~\ref{sec:preliminaries}.
The first main theorem of this paper is the following.

\begin{theorem}\label{thm:galois-pair-criterion}
Assume Conjecture~\ref{conj:wsc} \textup{(WSC)}.
Let \(K_{1}\) and \(K_{2}\) be number fields that are Galois over \(\QQ\).
Then the following two conditions are equivalent.
\begin{enumerate}[label=\textup{(\roman*)}]
\item \(r(K_{1}) > 0\), \(r(K_{2}) > 0\), and \(K_{1}^{+} \ne K_{2}^{+}\).
\item The regulators \(\reg(K_{1})\) and \(\reg(K_{2})\) are algebraically independent over \(\QQ\).
\end{enumerate}
\end{theorem}

In Theorem~\ref{thm:galois-pair-criterion}, we prove that \textup{(ii)} implies \textup{(i)} without assuming WSC.

The second main theorem of this paper is the following.

\begin{theorem}\label{thm:main-algebraic-independence}
Assume Conjecture~\ref{conj:wsc} \textup{(WSC)}.
Let \(m \ge 1\) be an integer, and let \(K_{1}, \dots, K_{m}\) be number fields that are Galois over \(\QQ\).
Suppose that
\[
r(K_{i}) > 0
\]
for each \(i \in \{1, \dots, m\}\), and assume that
\[
r(K_{i} \cap K_{j}) = 0
\]
for all distinct \(i, j \in \{1, \dots, m\}\).
Then the regulators
\[
\reg(K_{1}), \dots, \reg(K_{m})
\]
are algebraically independent over \(\QQ\).
\end{theorem}

We can rephrase the assumption on the intersections in Theorem~\ref{thm:main-algebraic-independence} in terms of maximal real subfields.
Indeed, as we show in Lemma~\ref{lem:intersection-rank-zero} of Section~\ref{sec:maximal-real-subfield}, for all distinct \(i, j\), the condition
\[
r(K_{i} \cap K_{j}) = 0
\]
is equivalent to
\[
K_{i}^{+} \cap K_{j}^{+} = \QQ
\]
and also to the condition that \(K_{i} \cap K_{j}\) is either \(\QQ\) or an imaginary quadratic field.

When \(m = 1\), Theorem~\ref{thm:main-algebraic-independence} asserts that
the regulator of a number field that is Galois over \(\QQ\) and whose unit group has positive rank is a transcendental number.
Under the weak Schanuel conjecture, the transcendence of the regulator in the case where the unit group has positive rank follows from the argument of Gun, Murty, and Rath~\cite[proof of Theorem~5.1]{GunMurtyRath2011}.
The case \(m = 1\) of Theorem~\ref{thm:main-algebraic-independence} reproduces this result for number fields that are Galois over \(\QQ\).
In this sense, the new content of Theorem~\ref{thm:main-algebraic-independence} lies in the algebraic independence in the case \(m \ge 2\).

Suppose that \(K_{1}, \dots, K_{m}\) in Theorem~\ref{thm:main-algebraic-independence} are all totally real. Then the intersection of any two of these fields is also totally real,
and hence the unit group of the intersection has rank \(0\) if and only if the intersection is \(\QQ\).
Therefore, in this case, Theorem~\ref{thm:main-algebraic-independence}
strengthens Cruz's linear independence for finite families~\cite[Theorem~5.3]{Cruz2025} to algebraic independence,
under the additional assumptions that each field is Galois over \(\QQ\) and that the intersection of any two distinct fields is \(\QQ\).
Indeed, since \(\QQalg/\QQ\) is an algebraic extension,
algebraic independence over \(\QQ\) is equivalent to algebraic independence over \(\QQalg\)~\cite[Theorem~4.3.5]{Roman2006},
and the latter implies linear independence over \(\QQalg\).
Note that, under the assumptions of Theorem~\ref{thm:main-algebraic-independence}, each \(K_{i}\) is different from \(\QQ\) because \(r(K_{i}) > 0\),
and \(K_{i} \ne K_{j}\) whenever \(i \ne j\) because \(r(K_{i} \cap K_{j}) = 0 < r(K_{i})\).
Since two distinct number fields that are Galois over \(\QQ\) are not isomorphic,
the result of Perlis~\cite{Perlis1977} shows that their Dedekind zeta functions are distinct, and thus the hypothesis of Cruz's theorem is satisfied.
In addition, Theorem~\ref{thm:main-algebraic-independence} does not assume that the degrees of \(K_{1}, \dots, K_{m}\) are equal to one another.
The fields \(K_{1}, \dots, K_{m}\) may include both totally real number fields and totally imaginary number fields.
We do not assume that the totally imaginary number fields are CM fields.

We characterize the transcendence of the ratio of two regulators as follows, including the case where the unit group of one or both of the fields has rank \(0\).

\begin{corollary}\label{cor:ratio-transcendence}
Assume Conjecture~\ref{conj:wsc} \textup{(WSC)}.
Let \(K_{1}\) and \(K_{2}\) be number fields that are Galois over \(\QQ\).
Then the following three conditions are equivalent.
\begin{enumerate}[label=\textup{(\roman*)}]
\item \(K_{1}^{+} \ne K_{2}^{+}\).
\item The ratio \(\reg(K_{1})/\reg(K_{2})\) is a transcendental number.
\item \(\reg(K_{1})/\reg(K_{2}) \notin \QQ\).
\end{enumerate}
In particular, the ratio \(\reg(K_{1})/\reg(K_{2})\) is either a rational number or a transcendental number, and it is a rational number if and only if \(K_{1}^{+} = K_{2}^{+}\).
\end{corollary}

In Corollary~\ref{cor:ratio-transcendence}, we prove that \textup{(ii)} implies \textup{(iii)} and that \textup{(iii)} implies \textup{(i)} without assuming WSC.

If \(K_{1}\) and \(K_{2}\) in Theorem~\ref{thm:galois-pair-criterion} are both totally real, then \(K_{i}^{+} = K_{i}\) for each \(i \in \{1, 2\}\), and we obtain the following corollary.

\begin{corollary}\label{cor:totally-real-pair}
Assume Conjecture~\ref{conj:wsc} \textup{(WSC)}.
Let \(K_{1}\) and \(K_{2}\) be totally real number fields different from \(\QQ\),
both of which are Galois over \(\QQ\).
Then the following two conditions are equivalent.
\begin{enumerate}[label=\textup{(\roman*)}]
\item \(K_{1} \ne K_{2}\).
\item The regulators \(\reg(K_{1})\) and \(\reg(K_{2})\) are algebraically independent over \(\QQ\).
\end{enumerate}
\end{corollary}

Corollary~\ref{cor:totally-real-pair} imposes no assumption on \(K_{1} \cap K_{2}\).
Hence, when two distinct totally real number fields that are Galois over \(\QQ\) are both different from \(\QQ\),
Corollary~\ref{cor:totally-real-pair} strengthens Cruz's linear independence over \(\QQalg\)~\cite[Theorem~5.3]{Cruz2025} to algebraic independence without any additional assumption on the intersection.

The proof of Theorem~\ref{thm:main-algebraic-independence} and the proof that \textup{(i)} implies \textup{(ii)} in Theorem~\ref{thm:galois-pair-criterion} follow a common strategy.
First, we express each regulator as the absolute value of the value of a nonzero homogeneous polynomial with rational coefficients at the logarithms of units that are positive real numbers.
Next, we show that these logarithms are linearly independent over \(\QQ\), and we deduce from WSC that they are algebraically independent over \(\QQ\).
Finally, we show by a comparison of multidegrees that these polynomials are algebraically independent over \(\QQ\).
This fact and the algebraic independence of the logarithms yield the algebraic independence of the regulators.
In Theorem~\ref{thm:main-algebraic-independence}, the polynomials corresponding to distinct fields have no indeterminates in common.
In Theorem~\ref{thm:galois-pair-criterion}, by contrast, if we put \(F \coloneqq K_{1} \cap K_{2}\) in the notation of that theorem, then \(r(F) > 0\) may occur.
For this reason, for each \(i \in \{1, 2\}\), we consider the units of \(F\) that are positive real numbers
and the units of \(K_{i}\) that are positive real numbers and whose norm relative to \(K_{i}/F\) is \(1\).
By means of the logarithms of these two kinds of units, we express \(\reg(K_{i})\) as the absolute value of the value of a multihomogeneous polynomial in two tuples of indeterminates.
The two polynomials then share the tuple of indeterminates corresponding to the units of \(F\).
Moreover, if \(K_{1}^{+} \ne K_{2}^{+}\), then at least one of \(r(K_{1}) - r(F)\) and \(r(K_{2}) - r(F)\) is positive.
The comparison of multidegrees uses this fact.
We give the definitions of homogeneous polynomials, multidegrees, and multihomogeneous polynomials in Section~\ref{sec:preliminaries}.

The organization of this paper is as follows.
In Section~\ref{sec:preliminaries}, we fix the notation and formulate the weak Schanuel conjecture.
In Section~\ref{sec:maximal-real-subfield}, we study the relation of the maximal real subfield of a number field to the rank of the unit group and to the regulator.
In particular, we treat the case where the maximal real subfields of two fields coincide, and we prove a lemma that we use in the proof that \textup{(ii)} implies \textup{(i)} in Theorem~\ref{thm:galois-pair-criterion}.
In Sections~\ref{sec:regulator-polynomial},~\ref{sec:log-linear-independence}, and~\ref{sec:polynomial-independence}, we prove, in this order, the lemmas corresponding to the three steps of the strategy described above.
In each of these sections, we first state the lemmas that we use in the proof of Theorem~\ref{thm:main-algebraic-independence} and then those that we use in the proof of Theorem~\ref{thm:galois-pair-criterion}.
In Section~\ref{sec:main}, we prove Theorem~\ref{thm:main-algebraic-independence} and Theorem~\ref{thm:galois-pair-criterion} in this order,
and we deduce Corollary~\ref{cor:ratio-transcendence} and Corollary~\ref{cor:totally-real-pair}.
Finally, in Section~\ref{sec:remarks}, we examine in turn the conditions on the Galois property, the rank of the unit group, the maximal real subfields, and the intersections.
We present examples in which the conclusions fail when we remove an assumption or weaken it in a specific way, as well as examples that illustrate Corollary~\ref{cor:ratio-transcendence}.

\section{Preliminaries}\label{sec:preliminaries}

In this section, for the convenience of the reader, we summarize the notation and terminology used in this paper and recall the basic facts needed in the later arguments.

We write \(\overline{z}\) for the complex conjugate of a complex number \(z\). The symbol \(\log\) denotes the natural logarithm of positive real numbers.

For a subset \(D\) of \(\RR\), we put
\[
\begin{aligned}
D_{ > 0}
& \coloneqq D \cap \{x \in \RR : x > 0\}, \\
D_{ \ge 0}
& \coloneqq D \cap \{x \in \RR : x \ge 0\}.
\end{aligned}
\]

For a finite extension \(L/K\) of fields, we write \([L:K]\) for its degree.
When an extension \(L/K\) of fields is a Galois extension, we write \(\Gal(L/K)\) for its Galois group.

Let \(\Omega\) be a field, let \(m \ge 1\) be an integer, and let \(K_{1}, \dots, K_{m}\) be subfields of \(\Omega\).
We write \(K_{1}\cdots K_{m}\) for the smallest subfield of \(\Omega\) that contains all of \(K_{1}, \dots, K_{m}\),
and we call it the compositum of \(K_{1}, \dots, K_{m}\).
In particular, we write \(KK'\) for the compositum of two fields \(K, K'\).

Let \(L/K\) be a finite extension of fields. For each \(x \in L\), consider the multiplication map by \(x\),
\[
\psi_{x} \colon L \longrightarrow L,
\qquad y \longmapsto xy.
\]
This map is \(K\)-linear. We define the norm map~\cite[Chapter~2, Review of norms and traces]{MilneANT} by
\[
\Norm_{L/K} \colon L \longrightarrow K,
\qquad x \longmapsto \det(\psi_{x}).
\]
The norm map is multiplicative: for all \(x, y \in L\), we have
\[
\Norm_{L/K}(xy) = \Norm_{L/K}(x)\Norm_{L/K}(y).
\]
Moreover, for every \(a \in K\), we have
\[
\Norm_{L/K}(a) = a^{[L:K]}.
\]
Furthermore, let \(K \subseteq L \subseteq E\) be fields such that \(E/K\) is a finite extension.
Then the norm is transitive; that is, for every \(x \in E\), we have
\[
\Norm_{E/K}(x) = \Norm_{L/K}\bigl(\Norm_{E/L}(x)\bigr).
\]
In addition, when \(L/K\) is a Galois extension, we have
\[
\Norm_{L/K}(x) = \prod_{\sigma \in \Gal(L/K)} \sigma(x)
\]
for every \(x \in L\).

The degree of a number field \(K\) means the degree \([K:\QQ]\) of \(K\) over \(\QQ\).
We say that a number field \(K\) is Galois over \(\QQ\) if \(K/\QQ\) is a Galois extension.
We call a number field \(K\) that is Galois over \(\QQ\) a cyclic field if \(\Gal(K/\QQ)\) is a cyclic group.
Furthermore, we call a quartic field \(K\) that is Galois over \(\QQ\) a biquadratic field if \(\Gal(K/\QQ)\) is isomorphic to the direct product of two cyclic groups of order \(2\).

We call a number field \(K\) a real number field if \(K \subseteq \RR\), and an imaginary number field otherwise.
For a number field \(K\), the field \(K \cap \RR\) is the largest real number field contained in \(K\). We call it the maximal real subfield of \(K\) and write
\[
K^{+} \coloneqq K \cap \RR.
\]
A real quadratic field means a quadratic field contained in \(\RR\), and an imaginary quadratic field means a quadratic field not contained in \(\RR\).

For a number field \(K\), we write \(\mathcal{O}_{K}\) for the ring of all algebraic integers contained in \(K\), that is, the ring of integers of \(K\);
we write \(\mathcal{O}_{K}^{\times}\) for its unit group and \(\mu(K)\) for the group of all roots of unity contained in \(K\).

For an extension \(L/K\) of number fields, the norm of an algebraic integer is an algebraic integer, and hence
\[
\Norm_{L/K}(\mathcal{O}_{L}) \subseteq \mathcal{O}_{K}.
\]
Moreover, the multiplicativity of the norm gives
\[
\Norm_{L/K}(\mathcal{O}_{L}^{\times}) \subseteq \mathcal{O}_{K}^{\times}
\]
and
\[
\Norm_{L/K}(\mu(L)) \subseteq \mu(K).
\]

In this paper, an embedding of a number field \(K\) means a field embedding \(\sigma \colon K \hookrightarrow \CC\).
We call \(\sigma\) a real embedding if \(\sigma(K) \subseteq \RR\), and a complex embedding otherwise.
When \(\sigma\) is a real embedding, we also write \(\sigma \colon K \hookrightarrow \RR\) to make explicit that \(\sigma(K) \subseteq \RR\).
For an embedding \(\sigma \colon K \hookrightarrow \CC\), we define the embedding \(\overline{\sigma} \colon K \hookrightarrow \CC\) by
\[
\overline{\sigma}(x) \coloneqq \overline{\sigma(x)} \quad (x \in K).
\]
We call \(\overline{\sigma}\) the complex conjugate of \(\sigma\).
An embedding \(\sigma\) is a real embedding if and only if \(\overline{\sigma} = \sigma\).
When \(\sigma\) is a complex embedding, \(\overline{\sigma}\) is also a complex embedding and \(\overline{\sigma} \ne \sigma\).
In this case, we call the unordered pair \(\{\sigma, \overline{\sigma}\}\) a complex conjugate pair of complex embeddings.

We write \(s(K)\) for the number of real embeddings of a number field \(K\) and \(t(K)\) for the number of complex conjugate pairs of complex embeddings of \(K\).
By~\cite[Chapter~4, Finiteness of the class number; Chapter~5, Statement of the theorem]{MilneANT}, we have
\begin{equation}
\label{eq:degree-embeddings} [K:\QQ] = s(K) + 2t(K).
\end{equation}
We call a number field totally real if all of its embeddings are real embeddings, and totally imaginary if all of its embeddings are complex embeddings.
We call a number field a CM field if it is a totally imaginary quadratic extension of a totally real number field.

Let \(K/F\) be an extension of number fields, and let \(\tau \colon F \hookrightarrow \CC\) be an embedding of \(F\).
We call an embedding \(\sigma\) of \(K\) an extension of \(\tau\) if \(\sigma|_{F} = \tau\).
Since \(K/F\) is a finite extension of fields of characteristic \(0\), it is a separable extension.
By the primitive element theorem, we can write \(K = F(\theta)\) for some \(\theta \in K\).
The minimal polynomial \(f\) of \(\theta\) over \(F\) has degree \([K:F]\) and is separable.
Therefore, the polynomial obtained by applying \(\tau\) to the coefficients of \(f\) has \([K:F]\) distinct roots in \(\CC\).
The extensions of \(\tau\) are in one-to-one correspondence with these roots~\cite[Proposition~2.4(b)]{MilneFT}.
Hence \(\tau\) has exactly \([K:F]\) extensions.

For a number field \(K\), the compositum of the images \(\sigma(K)\) of all embeddings \(\sigma \colon K \hookrightarrow \CC\) of \(K\)
is a number field that is Galois over \(\QQ\) and contains \(K\), and every number field that is Galois over \(\QQ\) and contains \(K\) contains this compositum~\cite[Chapter~V, \S~4, p.~242]{Lang2002}.
In other words, this compositum is the smallest number field that is Galois over \(\QQ\) and contains \(K\).
We call it the Galois closure of \(K\) over \(\QQ\).
By minimality, the field \(K\) uniquely determines its Galois closure.

Let \(G\) be an abelian group, written multiplicatively, with identity element \(1\).
Let \(n \ge 1\) be an integer. We say that \(g_{1}, \dots, g_{n} \in G\) are multiplicatively independent if,
for all \(A_{1}, \dots, A_{n} \in \ZZ\), the relation
\[
g_{1}^{A_{1}} \cdots g_{n}^{A_{n}} = 1
\]
implies
\[
A_{1} = \dots = A_{n} = 0.
\]

For a finitely generated abelian group \(G\), the quotient group of \(G\) by the subgroup consisting of all its elements of finite order is a finitely generated free abelian group.
We call the rank of this quotient group the rank of \(G\).
The rank is additive; that is, if \(G'\) is a subgroup of \(G\), then the rank of \(G\) equals the sum of the rank of \(G'\) and the rank of \(G/G'\).

For a number field \(K\), every root of unity contained in \(K\) is an algebraic integer, and its inverse is also a root of unity; hence \(\mu(K)\) is the subgroup of \(\mathcal{O}_{K}^{\times}\) consisting of all elements of finite order.
By Dirichlet's unit theorem~\cite[Theorem~5.1]{MilneANT}, the group \(\mathcal{O}_{K}^{\times}\) is a finitely generated abelian group, and therefore \(\mathcal{O}_{K}^{\times}/\mu(K)\) is a finitely generated free abelian group.
We write \(r(K)\) for its rank, that is, the rank of the unit group \(\mathcal{O}_{K}^{\times}\). We then have
\begin{equation}\label{eq:unit-rank}
r(K) = s(K) + t(K) - 1.
\end{equation}
We call a tuple of units whose cosets form a \(\ZZ\)-basis of \(\mathcal{O}_{K}^{\times}/\mu(K)\) a system of fundamental units of \(K\).
In particular, when \(r(K) = 1\), we call a unit whose coset forms a \(\ZZ\)-basis of \(\mathcal{O}_{K}^{\times}/\mu(K)\) a fundamental unit of \(K\).

Next, we define the logarithmic matrix and the regulator of a number field \(K\).
Put
\[
s \coloneqq s(K),
\qquad t \coloneqq t(K),
\qquad r \coloneqq r(K).
\]
Equation~\eqref{eq:unit-rank} gives \(r + 1 = s + t\).
Let \(\sigma_{1}, \dots, \sigma_{s}\) be the real embeddings of \(K\),
and let \(\sigma_{s + 1}, \dots, \sigma_{s + t}\) be embeddings obtained by choosing one from each complex conjugate pair of complex embeddings.
We call the tuple
\[
(\sigma_{1}, \dots, \sigma_{r + 1})
\]
obtained in this way a system of representatives of the embeddings of \(K\).
We define the weight of each \(\sigma_{a}\) by
\[
w_{a} \coloneqq
\begin{cases}
1, & 1 \le a \le s,\\
2, & s < a \le s + t.
\end{cases}
\]

When \(r > 0\), for \(u_{1}, \dots, u_{r} \in \mathcal{O}_{K}^{\times}\), we put
\[
\mathcal{M}_{K}(u_{1}, \dots, u_{r}) \coloneqq \bigl(w_{a}\log|\sigma_{a}(u_{j})|\bigr)_{1 \le a \le r + 1,\ 1 \le j \le r}.
\]
We call this matrix the logarithmic matrix of \(K\) associated with \(u_{1}, \dots, u_{r}\).
In addition, for \(a_{\ast} \in \{1, \dots, r + 1\}\), we write
\[
\mathcal{M}_{K}^{(a_{\ast})}(u_{1}, \dots, u_{r})
\]
for the \(r \times r\) matrix obtained from \(\mathcal{M}_{K}(u_{1}, \dots, u_{r})\) by removing the \(a_{\ast}\)-th row.
Let \(\varepsilon_{1}, \dots, \varepsilon_{r}\) be a system of fundamental units of \(K\), and put
\[
M_{K} \coloneqq \mathcal{M}_{K}(\varepsilon_{1}, \dots, \varepsilon_{r}).
\]
We call \(M_{K}\) a regulator matrix of \(K\).
We also choose \(a_{\ast} \in \{1, \dots, r + 1\}\), put
\[
M_{K}^{(a_{\ast})} \coloneqq \mathcal{M}_{K}^{(a_{\ast})}(\varepsilon_{1}, \dots, \varepsilon_{r}),
\]
and define
\[
\reg(K) \coloneqq |\det M_{K}^{(a_{\ast})}|.
\]

On the other hand, when \(r = 0\), we define
\[
\reg(K) \coloneqq 1.
\]

We call the number \(\reg(K)\) defined above the regulator of \(K\).

When \(r > 0\), for all \(u_{1}, \dots, u_{r} \in \mathcal{O}_{K}^{\times}\),
the value \(\bigl|\det\mathcal{M}_{K}^{(a_{\ast})}(u_{1}, \dots, u_{r})\bigr|\) does not depend on the choice of the index \(a_{\ast}\) of the removed row~\cite[\S~5.5, pp.~127--128]{Keune2023}.
In particular, \(\reg(K)\) does not depend on the choice of \(a_{\ast}\).
Furthermore, \(\reg(K)\) depends neither on the choice of the system of representatives of the embeddings nor on the choice of the system of fundamental units, and the field \(K\) uniquely determines \(\reg(K)\)~\cite[Appendix~A, \S~A.2, p.~184]{BiluBugeaudMignotte2014}.
We also have \(\reg(K) > 0\)~\cite[\S~5.4, p.~97]{Khanduja2022}.
In addition, if the cosets \(u_{1}\mu(K), \dots, u_{r}\mu(K)\) of \(u_{1}, \dots, u_{r} \in \mathcal{O}_{K}^{\times}\) are multiplicatively independent in \(\mathcal{O}_{K}^{\times}/\mu(K)\),
then the subgroup generated by these cosets has finite index, and \cite[Lemma~4.15]{Washington1997} gives
\begin{equation}
\label{eq:index-formula} \left|\det\mathcal{M}_{K}^{(a_{\ast})}(u_{1}, \dots, u_{r})\right| = \left[\mathcal{O}_{K}^{\times}/\mu(K):\left\langle u_{1}\mu(K), \dots, u_{r}\mu(K) \right\rangle\right]\reg(K)
\end{equation}
for every \(a_{\ast} \in \{1, \dots, r + 1\}\).
We call equation~\eqref{eq:index-formula} the index formula for logarithmic matrices.
The weight \(2\) for complex embeddings follows the usual convention for the regulator~\cite[Chapter~5]{MilneANT}.

\begin{remark}\label{rem:regulator-gram-determinant}
The definition of the regulator used in Cruz's paper~\cite[Definition~3.5]{Cruz2025} also gives the same value as \(\reg(K)\) defined above.
To verify this, we state the definition of that paper in the notation of the present paper.
In what follows, \(r\), \(s\), \(t\), the system of representatives \((\sigma_{1}, \dots, \sigma_{r + 1})\) of the embeddings,
and the weights \(w_{1}, \dots, w_{r + 1}\) are those used in the definition of the regulator.

We call the map
\[
\ell_{K} \colon \mathcal{O}_{K}^{\times} \longrightarrow \RR^{r + 1},
\qquad u \longmapsto \bigl(w_{a}\log|\sigma_{a}(u)|\bigr)_{1 \le a \le r + 1}
\]
the logarithmic embedding of \(K\),
and we call its image \(\Lambda_{K} \coloneqq \ell_{K}(\mathcal{O}_{K}^{\times})\) the log-unit lattice of \(K\)~\cite[Definition~3.4]{Cruz2025}.

That paper defines the Gram matrix of \(\Lambda_{K}\) with respect to a \(\ZZ\)-basis by means of the standard inner product on \(\RR^{r + 1}\),
and writes \(\det(\Lambda_{K})\) for its determinant. This value does not depend on the choice of the \(\ZZ\)-basis.
That paper then defines the regulator by
\[
\sqrt{\frac{\det(\Lambda_{K})}{s + t}}.
\]
Note, however, that the same paper writes \(r\) for the number of real embeddings and \(s\) for the number of complex conjugate pairs of complex embeddings,
and thus \(r + s\) there corresponds to \(s + t\) in the present paper.

The map \(\ell_{K}\) is a group homomorphism.
Moreover, since \(|\sigma_{a}(\zeta)| = 1\) for every \(\zeta \in \mu(K)\) and every \(a \in \{1, \dots, r + 1\}\),
the map \(\ell_{K}\) is trivial on \(\mu(K)\).

We first consider the case \(r > 0\).
Let \(\varepsilon_{1}, \dots, \varepsilon_{r}\) be the system of fundamental units used in the definition of the regulator.
Since \(\varepsilon_{1}, \dots, \varepsilon_{r}\) and \(\mu(K)\) generate \(\mathcal{O}_{K}^{\times}\),
the vectors \(\ell_{K}(\varepsilon_{1}), \dots, \ell_{K}(\varepsilon_{r})\),
that is, the columns of \(M_{K}\), generate \(\Lambda_{K}\).
In addition, \(\reg(K) > 0\) implies \(\det M_{K}^{(a_{\ast})} \ne 0\), and hence the columns of \(M_{K}\) are linearly independent over \(\RR\).
Therefore, the columns of \(M_{K}\) form a \(\ZZ\)-basis of \(\Lambda_{K}\). If we write \(M_{K}^{\mathsf{T}}\) for the transpose of \(M_{K}\),
then the Gram matrix with respect to this \(\ZZ\)-basis is \(M_{K}^{\mathsf{T}}M_{K}\).
As stated above, we have
\[
\bigl|\det\mathcal{M}_{K}^{(a)}(\varepsilon_{1}, \dots, \varepsilon_{r})\bigr| = \reg(K)
\]
for every \(a \in \{1, \dots, r + 1\}\), and thus the Cauchy--Binet formula~\cite[Proposition~3.4]{Serre2010} yields
\[
\begin{aligned}
\det(\Lambda_{K})
&= \det\bigl(M_{K}^{\mathsf{T}}M_{K}\bigr) \\
&= \sum_{a = 1}^{r + 1}\det\bigl(\mathcal{M}_{K}^{(a)}(\varepsilon_{1}, \dots, \varepsilon_{r})\bigr)^{2} \\
&= (r + 1)\reg(K)^{2}.
\end{aligned}
\]
Since \(r + 1 = s + t\) by equation~\eqref{eq:unit-rank} and \(\reg(K) > 0\), the value given by the definition in terms of the Gram matrix equals \(\reg(K)\).

When \(r = 0\), the group \(\mathcal{O}_{K}^{\times}/\mu(K)\) is a free abelian group of rank \(0\), that is, the trivial group, and hence \(\mathcal{O}_{K}^{\times} = \mu(K)\).
Since \(\ell_{K}\) is trivial on \(\mu(K)\), we have \(\Lambda_{K} = \{0\}\).
In this case, the \(\ZZ\)-basis of \(\Lambda_{K}\) is empty and \(s + t = 1\);
consequently, under the convention that the determinant of the empty Gram matrix is \(1\), the value given by the definition in terms of the Gram matrix is \(1\).
Since \(\reg(K) = 1\) by definition, this value equals \(\reg(K)\).
\end{remark}

Let \(R\) be a commutative ring.
Let \(d \ge 0\) be an integer, and let \(X_{1}, \dots, X_{d}\) be distinct indeterminates.
When \(d = 0\), there are no indeterminates, and the polynomial ring \(R[X_{1}, \dots, X_{d}]\) is understood to be \(R\) itself.
For a tuple \((n_{1}, \dots, n_{d})\) of nonnegative integers, we call the product
\[
\prod_{k = 1}^{d} X_{k}^{n_{k}}
\]
a monomial in \(X_{1}, \dots, X_{d}\), and we call the sum
\[
\sum_{k = 1}^{d} n_{k}
\]
the total degree of this monomial.
When \(d = 0\), the only tuple \((n_{1}, \dots, n_{d})\) of nonnegative integers is the empty tuple, the above product is the empty product, and the above sum is the empty sum.
Hence, in this case, the only monomial is \(1\), and the total degree of the monomial \(1\) is \(0\).
We call a nonzero polynomial \(P \in R[X_{1}, \dots, X_{d}]\) a homogeneous polynomial if all monomials appearing in \(P\) have the same total degree, and we call this common value the total degree of \(P\).
We regard the zero polynomial as a homogeneous polynomial having every nonnegative integer as its total degree.
We call a homogeneous polynomial of total degree \(1\) a homogeneous linear form.
That is, a homogeneous linear form is a polynomial of the form
\[
\sum_{k = 1}^{d} b_{k}X_{k}
\]
with \(b_{1}, \dots, b_{d} \in R\).

Let \(m \ge 1\) be an integer.
For each \(i \in \{1, \dots, m\}\), let \(d_{i} \ge 0\) be an integer, and let
\[
\boldsymbol{X}_{i} = (X_{i, 1}, \dots, X_{i, d_{i}})
\]
be a tuple of indeterminates.
If \(d_{i} = 0\), then we take \(\boldsymbol{X}_{i}\) to be the empty tuple.
We assume that, for all distinct \(i, j \in \{1, \dots, m\}\), the tuples \(\boldsymbol{X}_{i}\) and \(\boldsymbol{X}_{j}\) have no component in common.

Consider a monomial
\[
\prod_{i = 1}^{m}\left(\prod_{k = 1}^{d_{i}} X_{i, k}^{n_{i, k}}\right)
\]
in the indeterminates that are components of \(\boldsymbol{X}_{1}, \dots, \boldsymbol{X}_{m}\),
where each \(n_{i, k}\) is a nonnegative integer.
For each \(i \in \{1, \dots, m\}\), we call the sum
\[
\sum_{k = 1}^{d_{i}} n_{i, k}
\]
the total degree of this monomial with respect to \(\boldsymbol{X}_{i}\).
If \(d_{i} = 0\), then the above sum is the empty sum, and hence the total degree of every monomial with respect to \(\boldsymbol{X}_{i}\) is \(0\).
We also call the tuple of nonnegative integers
\[
\left( \sum_{k = 1}^{d_{1}} n_{1, k}, \dots, \sum_{k = 1}^{d_{m}} n_{m, k} \right)
\]
the multidegree of this monomial with respect to \((\boldsymbol{X}_{1}, \dots, \boldsymbol{X}_{m})\).

We write \(R[\boldsymbol{X}_{1}, \dots, \boldsymbol{X}_{m}]\) for the polynomial ring over \(R\) in the indeterminates that are components of \(\boldsymbol{X}_{1}, \dots, \boldsymbol{X}_{m}\).
We call a nonzero polynomial \(P \in R[\boldsymbol{X}_{1}, \dots, \boldsymbol{X}_{m}]\) a multihomogeneous polynomial with respect to \((\boldsymbol{X}_{1}, \dots, \boldsymbol{X}_{m})\)
if all monomials appearing in \(P\) have the same multidegree with respect to \((\boldsymbol{X}_{1}, \dots, \boldsymbol{X}_{m})\), and we call this common value the multidegree of \(P\) with respect to \((\boldsymbol{X}_{1}, \dots, \boldsymbol{X}_{m})\).
When \((\boldsymbol{X}_{1}, \dots, \boldsymbol{X}_{m})\) is clear from the context, we omit the reference to it and simply speak of the multidegree and of multihomogeneous polynomials.

Let \(L\) be a field, let \(A\) be a commutative \(L\)-algebra, and let \(n \ge 1\) be an integer.
We say that \(x_{1}, \dots, x_{n} \in A\) are algebraically independent over \(L\) if
\[
H(x_{1}, \dots, x_{n}) \ne 0
\]
for every nonzero polynomial \(H \in L[T_{1}, \dots, T_{n}]\) in indeterminates \(T_{1}, \dots, T_{n}\).

In this paper, we use only the following form of the weak Schanuel conjecture, which concerns the natural logarithms of algebraic numbers that are positive real numbers.

\begin{conjecture}[WSC]\label{conj:wsc}
Let \(n \ge 1\) be an integer, and let \(\alpha_{1}, \dots, \alpha_{n}\) be algebraic numbers that are positive real numbers.
If \(\log\alpha_{1}, \dots, \log\alpha_{n}\) are linearly independent over \(\QQ\), then these logarithms are algebraically independent over \(\QQ\).
\end{conjecture}

In this paper, WSC refers to Conjecture~\ref{conj:wsc}.
WSC is a special case of the weak Schanuel conjecture that Cruz uses~\cite[Conjecture~4.2]{Cruz2025}.
Hence the results of this paper also hold under that conjecture.

\section{Maximal real subfields}\label{sec:maximal-real-subfield}

In this section, we study the relation of the maximal real subfield of a number field to the rank of the unit group and to the regulator.
We first characterize the number fields whose unit groups have rank \(0\).

\begin{lemma}\label{lem:rank-zero}
Let \(K\) be a number field.
\begin{enumerate}[label=\textup{(\roman*)}]
\item \(r(K) = 0\) if and only if \(K = \QQ\) or \(K\) is an imaginary quadratic field.
\item If \(K\) is Galois over \(\QQ\), then \(r(K) = 0\) if and only if \(K^{+} = \QQ\).
\item If \(K\) is Galois over \(\QQ\) and \(r(K) > 0\), then \(r(K^{+}) > 0\).
\end{enumerate}
\end{lemma}

\begin{proof}
We first prove \textup{(i)}.
Suppose that \(r(K) = 0\).
Equation~\eqref{eq:unit-rank} gives \(s(K) + t(K) = 1\), and since \(s(K)\) and \(t(K)\) are nonnegative integers, either \(s(K) = 1\) and \(t(K) = 0\), or \(s(K) = 0\) and \(t(K) = 1\).
In the former case, equation~\eqref{eq:degree-embeddings} gives \([K:\QQ] = 1\), and hence \(K = \QQ\).
In the latter case, equation~\eqref{eq:degree-embeddings} gives \([K:\QQ] = 2\). Moreover, since \(K\) has no real embedding, the inclusion map \(K \hookrightarrow \CC\) is not a real embedding, and thus \(K \not\subseteq \RR\). Therefore \(K\) is an imaginary quadratic field.

Conversely, suppose that \(K = \QQ\) or \(K\) is an imaginary quadratic field.
If \(K = \QQ\), then equation~\eqref{eq:degree-embeddings} gives \(s(K) + 2t(K) = 1\), and hence \(s(K) = 1\) and \(t(K) = 0\).
If, on the other hand, \(K\) is an imaginary quadratic field, then the inclusion map \(K \hookrightarrow \CC\) is a complex embedding because \(K \not\subseteq \RR\), and thus \(t(K) \ge 1\); since equation~\eqref{eq:degree-embeddings} gives \(s(K) + 2t(K) = 2\), we obtain \(s(K) = 0\) and \(t(K) = 1\).
In either case, we have \(s(K) + t(K) = 1\), and equation~\eqref{eq:unit-rank} gives \(r(K) = 0\).

Next, we prove \textup{(ii)}.
Since \(K/\QQ\) is a Galois extension, the restriction of complex conjugation
\[
\kappa \colon K \longrightarrow K,
\qquad x \longmapsto \overline{x}
\]
is well-defined, and \(\kappa \in \Gal(K/\QQ)\).
The field \(K^{+} = K \cap \RR\) is the set of all elements of \(K\) fixed by \(\kappa\), that is, the fixed field of \(\langle\kappa\rangle\).
Since the order of \(\kappa\) is at most \(2\), the Galois correspondence gives \([K:K^{+}] \le 2\).

Suppose that \(r(K) = 0\).
By \textup{(i)}, either \(K = \QQ\) or \(K\) is an imaginary quadratic field.
If \(K = \QQ\), then \(K^{+} = \QQ\).
If \(K\) is an imaginary quadratic field, then \(K^{+} \ne K\) because \(K \not\subseteq \RR\); hence \([K:K^{+}] = 2 = [K:\QQ]\), and consequently \(K^{+} = \QQ\).

Conversely, suppose that \(K^{+} = \QQ\).
Then \([K:\QQ] \le 2\).
If \([K:\QQ] = 1\), then \(K = \QQ\).
If \([K:\QQ] = 2\), then \(K \ne K^{+}\) and thus \(K \not\subseteq \RR\); hence \(K\) is an imaginary quadratic field.
In either case, \textup{(i)} gives \(r(K) = 0\).

Finally, we prove \textup{(iii)}.
Suppose that \(K/\QQ\) is a Galois extension and that \(r(K) > 0\).
By \textup{(ii)}, we have \(K^{+} \ne \QQ\).
Moreover, \(K^{+}\) is a real number field and thus is not an imaginary quadratic field.
Therefore, \textup{(i)} applied to \(K^{+}\) yields \(r(K^{+}) > 0\).
\end{proof}

By Lemma~\ref{lem:rank-zero}~\textup{(i)} and the definition of the regulator, the regulator of \(\QQ\) and the regulator of every imaginary quadratic field are equal to \(1\).

Next, for the intersection of two number fields that are Galois over \(\QQ\), we rephrase the condition that its unit group has rank \(0\) in terms of maximal real subfields.

\begin{lemma}\label{lem:intersection-rank-zero}
Let \(K_{1}\) and \(K_{2}\) be number fields that are Galois over \(\QQ\).
Then the following three conditions are equivalent.
\begin{enumerate}[label=\textup{(\roman*)}]
\item \(r(K_{1} \cap K_{2}) = 0\).
\item \(K_{1}^{+} \cap K_{2}^{+} = \QQ\).
\item Either \(K_{1} \cap K_{2} = \QQ\), or \(K_{1} \cap K_{2}\) is an imaginary quadratic field.
\end{enumerate}
\end{lemma}

\begin{proof}
Put \(F \coloneqq K_{1} \cap K_{2}\).
Since \(K_{1}/\QQ\) and \(K_{2}/\QQ\) are Galois extensions, \(F/\QQ\) is also a Galois extension.
In addition, the definition of the maximal real subfield gives
\[
\begin{aligned}
F^{+}
&= (K_{1} \cap K_{2}) \cap \RR \\
&= (K_{1} \cap \RR) \cap (K_{2} \cap \RR) \\
&= K_{1}^{+} \cap K_{2}^{+}.
\end{aligned}
\]
Therefore, Lemma~\ref{lem:rank-zero}~\textup{(ii)} applied to \(F\) shows that \textup{(i)} and \textup{(ii)} are equivalent.
Likewise, Lemma~\ref{lem:rank-zero}~\textup{(i)} applied to \(F\) shows that \textup{(i)} and \textup{(iii)} are equivalent.
\end{proof}

Here we note that a number field \(K\) that is Galois over \(\QQ\) is either totally real or totally imaginary.
Indeed, by the normality of \(K/\QQ\), the image of every embedding \(K \hookrightarrow \CC\) coincides with \(K\).
If a real embedding exists, then \(K \subseteq \RR\), and all embeddings are real embeddings.
If no real embedding exists, then \(K\) is totally imaginary by definition.
Equations~\eqref{eq:degree-embeddings} and~\eqref{eq:unit-rank} give the rank of the unit group of \(K\) as
\[
r(K) =
\begin{cases}
[K:\QQ] - 1, & \text{if } K \text{ is totally real},\\[2pt]
\dfrac{[K:\QQ]}{2} - 1, & \text{if } K \text{ is totally imaginary}.
\end{cases}
\]
Hence the condition \(r(K) > 0\) is equivalent to the condition that the degree is at least \(2\) in the totally real case, and to the condition that the degree is at least \(4\) in the totally imaginary case.

Next, for an extension \(K/F\) of number fields, we compare \(r(F)\) with \(r(K)\) and show that, if the two ranks are equal, then the maximal real subfields of \(K\) and \(F\) coincide.

\begin{lemma}\label{lem:relative-rank}
Let \(F\) and \(K\) be number fields with \(F \subseteq K\).
Then
\[
r(F) \le r(K).
\]
Moreover, if \(r(K) = r(F)\), then \(K^{+} = F^{+}\).
\end{lemma}

\begin{proof}
Put \(n \coloneqq [K:F]\).
As stated in Section~\ref{sec:preliminaries}, every embedding of \(F\) has exactly \(n\) extensions.

Let \(\tau \colon F \hookrightarrow \RR\) be a real embedding of \(F\).
If \(\sigma\) is an extension of \(\tau\), then \(\overline{\sigma}|_{F} = \overline{\tau} = \tau\), and hence \(\overline{\sigma}\) is also an extension of \(\tau\).
Therefore, the complex embeddings among the extensions of \(\tau\) split into complex conjugate pairs.
Let \(s_{\tau}\) be the number of real embeddings of \(K\) that are extensions of \(\tau\), and let \(t_{\tau}\) be the number of complex conjugate pairs of complex embeddings of \(K\) that consist of extensions of \(\tau\).
Then the number of complex embeddings among the extensions of \(\tau\) is \(2t_{\tau}\).
Since \(\tau\) has \(n\) extensions in total, we obtain
\begin{equation}
\label{eq:real-place-degree} s_{\tau} + 2t_{\tau} = n.
\end{equation}
Since \(s_{\tau}\) and \(t_{\tau}\) are both nonnegative integers, equation~\eqref{eq:real-place-degree} and \(n \ge 1\) give
\begin{equation}
\label{eq:termwise-real} s_{\tau} + t_{\tau} \ge 1.
\end{equation}

Let \(\varphi \colon F \hookrightarrow \CC\) be a complex embedding of \(F\).
Every extension \(\sigma\) of \(\varphi\) satisfies \(\varphi(F) \subseteq \sigma(K)\).
Since \(\varphi(F) \not\subseteq \RR\), we see that \(\sigma(K) \not\subseteq \RR\).
Hence \(\sigma\) is a complex embedding.
Each of \(\varphi\) and \(\overline{\varphi}\) has \(n\) extensions.
The map \(\sigma \longmapsto \overline{\sigma}\) is a bijection from the set of all extensions of \(\varphi\) onto the set of all extensions of \(\overline{\varphi}\).
Since \(\varphi \ne \overline{\varphi}\),
for the complex conjugate pair \(\{\varphi, \overline{\varphi}\}\) of complex embeddings of \(F\),
there are exactly \(n\) complex conjugate pairs of complex embeddings of \(K\) that consist of an extension of \(\varphi\) and an extension of \(\overline{\varphi}\).

Every embedding of \(K\) determines a unique embedding of \(F\) as its restriction to \(F\).
Moreover, the restriction to \(F\) of a real embedding of \(K\) is a real embedding. Therefore
\[
\begin{aligned}
s(K)
& = \sum_{\tau \colon F \hookrightarrow \RR} s_{\tau}, \\
t(K)
& = \left(\sum_{\tau \colon F \hookrightarrow \RR} t_{\tau}\right) + nt(F),
\end{aligned}
\]
where the sums run over all real embeddings \(\tau\) of \(F\).
Consequently,
\begin{equation}
\label{eq:places-count} s(K) + t(K) = \left(\sum_{\tau \colon F \hookrightarrow \RR}(s_{\tau} + t_{\tau})\right) + nt(F).
\end{equation}
In addition, since
\[
s(F) = \sum_{\tau \colon F \hookrightarrow \RR} 1,
\]
equations~\eqref{eq:unit-rank} and~\eqref{eq:places-count} give
\begin{equation}
\label{eq:relative-rank-difference}
\begin{aligned}
r(K) - r(F)
& = \bigl(s(K) + t(K)\bigr) - \bigl(s(F) + t(F)\bigr) \\
& = \left(\sum_{\tau \colon F \hookrightarrow \RR}(s_{\tau} + t_{\tau} - 1)\right) + (n - 1)t(F).
\end{aligned}
\end{equation}
By inequality~\eqref{eq:termwise-real}, we have
\[
s_{\tau} + t_{\tau} - 1 \ge 0
\]
for every real embedding \(\tau\) of \(F\).
Furthermore, \(n \ge 1\) and \(t(F) \ge 0\) give
\[
(n - 1)t(F) \ge 0.
\]
Therefore,
\[
r(F) \le r(K).
\]

Assume that \(r(K) = r(F)\).
Then the right-hand side of equation~\eqref{eq:relative-rank-difference} is a sum of nonnegative terms whose value is \(0\), and hence each term is \(0\).
That is, for each real embedding \(\tau\) of \(F\), we have
\begin{equation}
\label{eq:equality-real} s_{\tau} + t_{\tau} = 1,
\end{equation}
and \((n - 1)t(F) = 0\).

If \(t(F) > 0\), then \((n - 1)t(F) = 0\) gives \(n = 1\), and thus \(K = F\).
Hence \(K^{+} = F^{+}\).

In what follows, suppose that \(t(F) = 0\), that is, that \(F\) is totally real.
Since equation~\eqref{eq:degree-embeddings} gives \(s(F) = [F:\QQ] \ge 1\), there exists a real embedding \(\tau\) of \(F\).
Equations~\eqref{eq:real-place-degree} and~\eqref{eq:equality-real} give
\[
t_{\tau} = n - 1,
\qquad s_{\tau} = 2 - n,
\]
and \(s_{\tau} \ge 0\) implies \(n = 1\) or \(n = 2\).

If \(n = 1\), then \(K = F\), and therefore \(K^{+} = F^{+}\).

If \(n = 2\), then
\[
s_{\tau} = 0,
\qquad t_{\tau} = 1
\]
for each real embedding \(\tau\) of \(F\).
Since \(F\) is totally real, the expression for \(s(K)\) obtained above gives
\[
s(K) = \sum_{\tau \colon F \hookrightarrow \RR} s_{\tau} = 0,
\]
and thus \(K\) is totally imaginary.
Therefore \(K\) is a totally imaginary quadratic extension of the totally real number field \(F\), that is, a CM field.
Since \(F\) is totally real, the inclusion map \(F \hookrightarrow \CC\) is a real embedding, and \(F \subseteq \RR\).
Hence \(F^{+} = F\), and we have
\[
F \subseteq K \cap \RR = K^{+} \subseteq K.
\]
Since \([K:K^{+}][K^{+}:F] = [K:F] = 2\), either \(K^{+} = F\) or \(K^{+} = K\).
Since \(K\) is totally imaginary, the inclusion map \(K \hookrightarrow \CC\) is not a real embedding, and \(K \not\subseteq \RR\).
Thus \(K^{+} \ne K\).
Consequently,
\[
K^{+} = F = F^{+}.
\]

We conclude that, if \(r(K) = r(F)\), then \(K^{+} = F^{+}\).
\end{proof}

Furthermore, we introduce the group of all units that are positive real numbers, and we show that its rank equals the rank of the unit group of the maximal real subfield.
For a number field \(K\), we write \(U(K)\) for the set of all units of \(K\) that are positive real numbers; that is,
\[
U(K) \coloneqq \mathcal{O}_{K}^{\times} \cap \RR_{ > 0}.
\]
Since \(U(K)\) is a subgroup of the finitely generated abelian group \(\mathcal{O}_{K}^{\times}\), it is finitely generated.
Moreover, the multiplicative group of positive real numbers has no element of finite order other than \(1\),
and hence, by the structure theorem for finitely generated abelian groups, \(U(K)\) is a free abelian group.

Since \(K^{+} \subseteq \RR\), the sign map
\[
\mathcal{O}_{K^{+}}^{\times} \longrightarrow \{\pm 1\},
\qquad u \longmapsto \frac{u}{|u|}
\]
is well-defined, and it is a group homomorphism.
This map is surjective because \(-1 \in \mathcal{O}_{K^{+}}^{\times}\), and its kernel is \(U(K^{+})\).
Therefore, \(U(K^{+})\) is a subgroup of index \(2\) in \(\mathcal{O}_{K^{+}}^{\times}\),
and in particular the rank of \(\mathcal{O}_{K^{+}}^{\times}/U(K^{+})\) is \(0\).
By the additivity of the rank, the rank of \(U(K^{+})\) equals \(r(K^{+})\).

In addition, since \(K/K^{+}\) is an extension of number fields, we have
\[
\mathcal{O}_{K^{+}}^{\times} = \mathcal{O}_{K}^{\times} \cap K^{+},
\]
and the definition of \(K^{+}\) further yields
\[
\begin{aligned}
U(K)
&= \mathcal{O}_{K}^{\times} \cap \RR_{ > 0} \\
&= (\mathcal{O}_{K}^{\times} \cap K^{+}) \cap \RR_{ > 0} \\
&= \mathcal{O}_{K^{+}}^{\times} \cap \RR_{ > 0} \\
&= U(K^{+}).
\end{aligned}
\]
Thus \(U(K)\) is a free abelian group of rank \(r(K^{+})\).

Finally, we show that, when the maximal real subfields of two number fields that are Galois over \(\QQ\) coincide, the ratio of their regulators is a positive rational number.
For this purpose, we first study the relation between a number field that is Galois over \(\QQ\) and its maximal real subfield.

\begin{lemma}\label{lem:real-part-structure}
Let \(K\) be a number field that is Galois over \(\QQ\).
Then the following hold.
\begin{enumerate}[label=\textup{(\roman*)}]
\item If \(K^{+}/\QQ\) is not a Galois extension, then \(K\) is the Galois closure of \(K^{+}\) over \(\QQ\).
\item If \(K^{+}/\QQ\) is a Galois extension, then either \(K = K^{+}\),
or \(K\) is a totally imaginary quadratic extension of the totally real number field \(K^{+}\), that is, a CM field.
\end{enumerate}
\end{lemma}

\begin{proof}
Since \(K/\QQ\) is a Galois extension, the restriction of complex conjugation
\[
\kappa \colon K \longrightarrow K,
\qquad x \longmapsto \overline{x}
\]
is well-defined, and \(\kappa \in \Gal(K/\QQ)\). Moreover, \(K^{+}\) is the fixed field of \(\langle\kappa\rangle\); that is,
\[
K^{+} = K^{\langle\kappa\rangle} \coloneqq \left\{x \in K : \kappa(x) = x\right\}.
\]

If \(K \subseteq \RR\), then \(\kappa\) is the identity map of \(K\), and \(K^{+} = K\).
In this case, \(K^{+}/\QQ\) is a Galois extension; hence the assumption of \textup{(i)} is not satisfied, and \textup{(ii)} holds because \(K = K^{+}\).

In what follows, suppose that \(K \not\subseteq \RR\).
Then \(\kappa\) has order \(2\).
Since \(K^{+}\) is the fixed field of \(\langle\kappa\rangle\), the Galois correspondence gives
\[
[K:K^{+}] = 2.
\]

We first prove \textup{(i)}.
Suppose that \(K^{+}/\QQ\) is not a Galois extension, and let \(L\) be the Galois closure of \(K^{+}\) over \(\QQ\).
Since \(K\) is a number field that is Galois over \(\QQ\) and contains \(K^{+}\), the minimality of the Galois closure gives
\[
K^{+} \subseteq L \subseteq K.
\]
Hence
\[
[K:L][L:K^{+}] = [K:K^{+}] = 2.
\]
It follows that \([K:L] = 1\) or \([L:K^{+}] = 1\).
Therefore, \(L = K\) or \(L = K^{+}\).
Since \(L/\QQ\) is a Galois extension and \(K^{+}/\QQ\) is not, we see that \(L \ne K^{+}\).
Thus \(L = K\); that is, \(K\) is the Galois closure of \(K^{+}\) over \(\QQ\).

Next, we prove \textup{(ii)}.
Suppose that \(K^{+}/\QQ\) is a Galois extension.
Since \(K^{+}\) is a Galois extension of \(\QQ\) contained in \(\RR\), it is totally real.
Moreover, since \(K/\QQ\) is a Galois extension, \(K\) is either totally real or totally imaginary, and since \(K \not\subseteq \RR\), the field \(K\) is totally imaginary.
Furthermore, as shown above, \([K:K^{+}] = 2\).
Therefore, \(K\) is a totally imaginary quadratic extension of the totally real number field \(K^{+}\), that is, a CM field.
\end{proof}

\begin{lemma}\label{lem:same-real-part-ratio}
Let \(K_{1}\) and \(K_{2}\) be number fields that are Galois over \(\QQ\), and suppose that \(K_{1}^{+} = K_{2}^{+}\).
Then
\[
\frac{\reg(K_{1})}{\reg(K_{2})} \in \QQ_{ > 0}.
\]
\end{lemma}

\begin{proof}
Since \(K_{1}^{+} = K_{2}^{+}\), the extension \(K_{1}^{+}/\QQ\) is a Galois extension if and only if \(K_{2}^{+}/\QQ\) is a Galois extension.

Suppose first that \(K_{1}^{+}/\QQ\) is not a Galois extension. Then Lemma~\ref{lem:real-part-structure}~\textup{(i)} applied to \(K_{1}\) and \(K_{2}\) shows that
both \(K_{1}\) and \(K_{2}\) are the Galois closure over \(\QQ\) of the common maximal real subfield \(K_{1}^{+} = K_{2}^{+}\).
Since the Galois closure is unique, we have \(K_{1} = K_{2}\), and
\[
\frac{\reg(K_{1})}{\reg(K_{2})} = 1 \in \QQ_{ > 0}.
\]
Suppose next that \(K_{1}^{+}/\QQ\) is a Galois extension. We first show that
\[
\frac{\reg(K_{i})}{\reg(K_{i}^{+})} \in \QQ_{ > 0}
\]
for each \(i \in \{1, 2\}\).
Take an arbitrary \(i \in \{1, 2\}\).
If \(K_{i} = K_{i}^{+}\), then this ratio is \(1\).
If \(K_{i} \ne K_{i}^{+}\), then, by Lemma~\ref{lem:real-part-structure}~\textup{(ii)}, \(K_{i}\) is a totally imaginary quadratic extension of the totally real number field \(K_{i}^{+}\), that is, a CM field.
In this situation, if \(r(K_{i}) = 0\), then Lemma~\ref{lem:rank-zero}~\textup{(ii)} gives \(K_{i}^{+} = \QQ\); hence \(\reg(K_{i}) = \reg(K_{i}^{+}) = 1\) by definition, and this ratio is \(1\).
If \(r(K_{i}) > 0\), then this ratio equals \(2^{r(K_{i})}\) or \(2^{r(K_{i}) - 1}\)~\cite[Theorem~4.12 and Proposition~4.16]{Washington1997}.
In each case, this ratio is a positive rational number.
Since \(K_{1}^{+} = K_{2}^{+}\), we have \(\reg(K_{1}^{+}) = \reg(K_{2}^{+})\), and we conclude that
\[
\frac{\reg(K_{1})}{\reg(K_{2})} = \frac{\reg(K_{1})/\reg(K_{1}^{+})}{\reg(K_{2})/\reg(K_{2}^{+})} \in \QQ_{ > 0}.
\]
\end{proof}

\section{Polynomial expressions for regulators}\label{sec:regulator-polynomial}

In this section, we express the regulator of a number field that is Galois over \(\QQ\) and whose unit group has positive rank as the absolute value of the value of a nonzero homogeneous polynomial with rational coefficients at the logarithms of units that are positive real numbers.
Furthermore, when two number fields \(F \subseteq K\) that are Galois over \(\QQ\) satisfy \(r(K) > 0\), we express the regulator of \(K\) as the absolute value of the value of a nonzero multihomogeneous polynomial with rational coefficients at the logarithms of units that are positive real numbers.

\begin{lemma}\label{lem:polynomial}
Let \(K\) be a number field that is Galois over \(\QQ\), and suppose that \(r(K) > 0\).
Put \(r \coloneqq r(K)\) and \(r_{+} \coloneqq r(K^{+})\).
Lemma~\ref{lem:rank-zero}~\textup{(iii)} gives \(r_{+} > 0\).
Let \(\eta_{1}, \dots, \eta_{r_{+}}\) be a \(\ZZ\)-basis of \(U(K)\).
Put
\[
x_{k} \coloneqq \log\eta_{k} \quad (1 \le k \le r_{+}),
\]
and let \(X_{1}, \dots, X_{r_{+}}\) be indeterminates.
Then there exists a nonzero homogeneous polynomial
\[
P \in \ZZ[X_{1}, \dots, X_{r_{+}}]
\]
of total degree \(r\) with integer coefficients such that
\[
\reg(K) = |P(x_{1}, \dots, x_{r_{+}})|.
\]
\end{lemma}

\begin{proof}
Let \(\varepsilon_{1}, \dots, \varepsilon_{r}\) be a system of fundamental units of \(K\),
and let \((\sigma_{1}, \dots, \sigma_{r + 1})\) and \(w_{1}, \dots, w_{r + 1}\) be, respectively,
the system of representatives of the embeddings of \(K\) and the weights in the definition of the regulator.

For each \(a \in \{1, \dots, r + 1\}\) and each \(j \in \{1, \dots, r\}\), put
\[
\alpha_{a, j} \coloneqq |\sigma_{a}(\varepsilon_{j})|^{w_{a}}.
\]
Since \(K/\QQ\) is a Galois extension, \(\sigma_{a}(\varepsilon_{j})\) is a unit of \(K\).

If \(\sigma_{a}\) is a real embedding, then \(w_{a} = 1\) and
\[
\alpha_{a, j} = |\sigma_{a}(\varepsilon_{j})| = \pm\sigma_{a}(\varepsilon_{j}).
\]
The number \(|\sigma_{a}(\varepsilon_{j})|\) is a positive real number.
Moreover, both of \(\pm\sigma_{a}(\varepsilon_{j})\) are units of \(K\).
Therefore, \(\alpha_{a, j} \in U(K)\).

If \(\sigma_{a}\) is a complex embedding, then \(w_{a} = 2\) and
\[
\begin{aligned}
\alpha_{a, j}
& = |\sigma_{a}(\varepsilon_{j})|^{2}\\
& = \sigma_{a}(\varepsilon_{j})\overline{\sigma_{a}(\varepsilon_{j})}.
\end{aligned}
\]
The number \(|\sigma_{a}(\varepsilon_{j})|^{2}\) is a positive real number.
In addition, by the normality of \(K/\QQ\), complex conjugation induces an automorphism of \(K\),
and hence \(\overline{\sigma_{a}(\varepsilon_{j})}\) is a unit of \(K\).
Thus \(\sigma_{a}(\varepsilon_{j})\overline{\sigma_{a}(\varepsilon_{j})}\) is a unit of \(K\).
Therefore, \(\alpha_{a, j} \in U(K)\) in this case as well.

Consequently, \(\alpha_{a, j} \in U(K)\) in either case.

Since \(r > 0\), by definition, the \((a, j)\)-entry of the regulator matrix \(M_{K}\) is
\[
w_{a}\log|\sigma_{a}(\varepsilon_{j})| = \log\alpha_{a, j}.
\]

Since \(\eta_{1}, \dots, \eta_{r_{+}}\) is a \(\ZZ\)-basis of \(U(K)\),
for each \(a \in \{1, \dots, r + 1\}\) and each \(j \in \{1, \dots, r\}\),
there exist integers \(b_{a, j, k}\) \((1 \le k \le r_{+})\) such that
\[
\alpha_{a, j} = \prod_{k = 1}^{r_{+}}\eta_{k}^{b_{a, j, k}}.
\]
Hence
\[
\log\alpha_{a, j} = \sum_{k = 1}^{r_{+}}b_{a, j, k}x_{k}.
\]
We now define
\[
P(X_{1}, \dots, X_{r_{+}}) \coloneqq \det\left(\sum_{k = 1}^{r_{+}}b_{a, j, k}X_{k}\right)_{1 \le a, j \le r}.
\]
If we choose \(a_{\ast} = r + 1\) as the index of the row to remove in the definition of the regulator, then we obtain
\[
\begin{aligned}
\reg(K)
& = \left|\det\bigl(w_{a}\log|\sigma_{a}(\varepsilon_{j})|\bigr)_{1 \le a, j \le r}\right|\\
& = |P(x_{1}, \dots, x_{r_{+}})|.
\end{aligned}
\]
The polynomial \(P\) has integer coefficients, and each entry of the matrix used in the definition of \(P\) is a homogeneous linear form; hence \(P\) is a homogeneous polynomial of total degree \(r\).
Since \(\reg(K) > 0\), the polynomial \(P\) is not the zero polynomial.
\end{proof}

Next, we prepare the notation that we use for the expression by multihomogeneous polynomials.
Let \(K/F\) be an extension of number fields. We write \(U^{1}(K/F)\) for the set of all elements of \(U(K)\) whose norm relative to \(K/F\) is \(1\); that is,
\[
U^{1}(K/F) \coloneqq \left\{u \in U(K) : \Norm_{K/F}(u) = 1\right\}.
\]
By the multiplicativity of the norm, \(U^{1}(K/F)\) is a subgroup of \(U(K)\).
Hence, as a subgroup of the free abelian group \(U(K)\), the group \(U^{1}(K/F)\) is also a free abelian group.
Furthermore, we put
\[
\mathcal{E}(K/F) \coloneqq \ker\bigl(\Norm_{K/F} \colon \mathcal{O}_{K}^{\times} \longrightarrow \mathcal{O}_{F}^{\times}\bigr).
\]
That is, \(\mathcal{E}(K/F)\) is the group of all units of \(K\) whose norm relative to \(K/F\) is \(1\).
Unlike in the case of \(U^{1}(K/F)\), we do not require the elements of \(\mathcal{E}(K/F)\) to be positive real numbers.

\begin{lemma}\label{lem:relative-unit-kernel}
Let \(F\) and \(K\) be number fields with \(F \subseteq K\) and \(r(K) > 0\).
Put
\[
\delta \coloneqq r(F),
\qquad \rho \coloneqq r(K) - r(F).
\]
Then the rank of \(\mathcal{E}(K/F)\) is \(\rho\). In particular, \(\rho \ge 0\).

Let \(\varepsilon_{1}, \dots, \varepsilon_{\delta}\) be a system of fundamental units of \(F\),
and choose \(v_{1}, \dots, v_{\rho} \in \mathcal{E}(K/F)\) such that their cosets form a \(\ZZ\)-basis of
\[
\mathcal{E}(K/F)/\bigl(\mathcal{E}(K/F) \cap \mu(K)\bigr).
\]
For each of these two tuples, if the corresponding rank is \(0\), then we take the tuple to be the empty tuple.
Then the cosets
\[
\varepsilon_{1}\mu(K), \dots, \varepsilon_{\delta}\mu(K), v_{1}\mu(K), \dots, v_{\rho}\mu(K)
\]
are multiplicatively independent in \(\mathcal{O}_{K}^{\times}/\mu(K)\),
and the subgroup generated by these cosets has finite index.
\end{lemma}

\begin{proof}
For every \(u \in \mathcal{O}_{F}^{\times}\), we have
\[
\Norm_{K/F}(u) = u^{[K:F]}.
\]
Hence
\[
\left\{u^{[K:F]} : u \in \mathcal{O}_{F}^{\times}\right\} \subseteq \Norm_{K/F}\bigl(\mathcal{O}_{K}^{\times}\bigr) \subseteq \mathcal{O}_{F}^{\times}.
\]
Since \(\mathcal{O}_{F}^{\times}\) is a finitely generated abelian group, the leftmost group is a subgroup of finite index in \(\mathcal{O}_{F}^{\times}\).
Thus the groups at both ends have rank \(r(F) = \delta\).
By the additivity of the rank, \(\Norm_{K/F}\bigl(\mathcal{O}_{K}^{\times}\bigr)\) therefore also has rank \(\delta\).
Since \(\mathcal{E}(K/F)\) is the kernel of the group homomorphism \(\Norm_{K/F} \colon \mathcal{O}_{K}^{\times} \longrightarrow \mathcal{O}_{F}^{\times}\),
the additivity of the rank again shows that the rank of \(\mathcal{E}(K/F)\) is
\[
r(K) - \delta = \rho.
\]

Since \(\mathcal{E}(K/F)\) is a subgroup of the finitely generated abelian group \(\mathcal{O}_{K}^{\times}\), it is finitely generated.
Moreover, since \(\mu(K)\) is the subgroup of \(\mathcal{O}_{K}^{\times}\) consisting of all elements of finite order,
\(\mathcal{E}(K/F) \cap \mu(K)\) is the subgroup of \(\mathcal{E}(K/F)\) consisting of all elements of finite order.
Therefore, the quotient group \(\mathcal{E}(K/F)/\bigl(\mathcal{E}(K/F) \cap \mu(K)\bigr)\) is a finitely generated free abelian group,
and its rank equals the rank \(\rho\) of \(\mathcal{E}(K/F)\).
Hence we can choose \(v_{1}, \dots, v_{\rho}\) as in the statement.

Since \(\delta + \rho = r(K) > 0\), the tuple
\(\varepsilon_{1}\mu(K), \dots, \varepsilon_{\delta}\mu(K), v_{1}\mu(K), \dots, v_{\rho}\mu(K)\)
is nonempty.
We show that this tuple is multiplicatively independent.
Take arbitrary \(A_{1}, \dots, A_{\delta}, B_{1}, \dots, B_{\rho} \in \ZZ\), and assume that
\[
\prod_{j = 1}^{\delta}\varepsilon_{j}^{A_{j}}\prod_{j = 1}^{\rho}v_{j}^{B_{j}} \in \mu(K).
\]
We apply \(\Norm_{K/F}\) to this product.
Since \(\Norm_{K/F}(\mu(K)) \subseteq \mu(F)\) and
\(\Norm_{K/F}(v_{j}) = 1\) for each \(j \in \{1, \dots, \rho\}\), we obtain
\[
\left(\prod_{j = 1}^{\delta}\varepsilon_{j}^{A_{j}}\right)^{[K:F]} \in \mu(F).
\]
Since \(\mathcal{O}_{F}^{\times}/\mu(F)\) has no element of finite order other than the identity element, we obtain
\[
\prod_{j = 1}^{\delta}\varepsilon_{j}^{A_{j}} \in \mu(F).
\]
Since \(\varepsilon_{1}, \dots, \varepsilon_{\delta}\) is a system of fundamental units,
all \(A_{j}\) are \(0\).
Furthermore, by the choice of \(v_{1}, \dots, v_{\rho}\),
all \(B_{j}\) are also \(0\).

In addition, since
\[
\delta + \rho = r(K),
\]
the subgroup generated by these cosets has finite index in
\(\mathcal{O}_{K}^{\times}/\mu(K)\).
\end{proof}

\begin{lemma}\label{lem:relative-log-decomposition}
Let \(F\) and \(K\) be number fields that are Galois over \(\QQ\),
with \(F \subseteq K\) and \(r(K) > 0\).
Put
\[
\delta \coloneqq r(F),
\qquad \rho \coloneqq r(K) - r(F),
\qquad \delta_{+} \coloneqq r(F^{+}),
\]
and let \(e\) be the rank of \(U^{1}(K/F)\).
Let \(\eta_{1}, \dots, \eta_{\delta_{+}}\) be a \(\ZZ\)-basis of \(U(F)\),
let \(\xi_{1}, \dots, \xi_{e}\) be a \(\ZZ\)-basis of \(U^{1}(K/F)\), and put
\[
x_{k} \coloneqq \log\eta_{k} \quad (1 \le k \le \delta_{+}),
\qquad y_{k} \coloneqq \log\xi_{k} \quad (1 \le k \le e).
\]
Moreover, choose \(\varepsilon_{1}, \dots, \varepsilon_{\delta}\) and \(v_{1}, \dots, v_{\rho}\) as in Lemma~\ref{lem:relative-unit-kernel}.
For each of the above four tuples, if the corresponding rank is \(0\), then we take the tuple to be the empty tuple.

Let \((\sigma_{1}, \dots, \sigma_{r(K) + 1})\) be the system of representatives of the embeddings of \(K\) in the definition of the logarithmic matrix in Section~\ref{sec:preliminaries},
and let \(w_{1}, \dots, w_{r(K) + 1}\) be the corresponding weights.
Then there exist integers
\[
b_{a,j,k} \quad (1 \le a \le r(K) + 1,\ 1 \le j \le \delta,\ 1 \le k \le \delta_{+})
\]
and integers
\[
c_{a,j,k} \quad (1 \le a \le r(K) + 1,\ 1 \le j \le \rho,\ 1 \le k \le e)
\]
such that
\[
2w_{a}\log|\sigma_{a}(\varepsilon_{j})| = \sum_{k = 1}^{\delta_{+}}b_{a,j,k}x_{k} \quad (1 \le a \le r(K) + 1,\ 1 \le j \le \delta)
\]
and
\[
2w_{a}\log|\sigma_{a}(v_{j})| = \sum_{k = 1}^{e}c_{a,j,k}y_{k} \quad (1 \le a \le r(K) + 1,\ 1 \le j \le \rho).
\]
\end{lemma}

\begin{proof}
Since \(K/\QQ\) is a Galois extension, we can regard each \(\sigma_{a}\) as an element of \(\Gal(K/\QQ)\).
In addition, the restriction of complex conjugation
\[
\kappa \colon K \longrightarrow K,
\qquad x \longmapsto \overline{x}
\]
is well-defined, and \(\kappa \in \Gal(K/\QQ)\).
Since \(F/\QQ\) is also a Galois extension, we have \(\sigma_{a}(F) = F\) and \(\kappa(F) = F\).

For each \(a \in \{1, \dots, r(K) + 1\}\) and each \(j \in \{1, \dots, \delta\}\), we have
\[
|\sigma_{a}(\varepsilon_{j})|^{2w_{a}} = \bigl(\sigma_{a}(\varepsilon_{j})\kappa(\sigma_{a}(\varepsilon_{j}))\bigr)^{w_{a}}.
\]
Since \(\varepsilon_{j} \in \mathcal{O}_{F}^{\times}\) and \(F\) is stable under both \(\sigma_{a}\) and \(\kappa\),
both \(\sigma_{a}(\varepsilon_{j})\) and \(\kappa(\sigma_{a}(\varepsilon_{j}))\) belong to \(\mathcal{O}_{F}^{\times}\).
Furthermore, the above product is a positive real number, and hence
\[
|\sigma_{a}(\varepsilon_{j})|^{2w_{a}} \in U(F).
\]
Therefore, there exist integers \(b_{a,j,k}\) \((1 \le k \le \delta_{+})\) such that
\[
2w_{a}\log|\sigma_{a}(\varepsilon_{j})| = \sum_{k = 1}^{\delta_{+}}b_{a,j,k}x_{k}.
\]

Since \(F/\QQ\) is a Galois extension, \(\Gal(K/F)\) is a normal subgroup of \(\Gal(K/\QQ)\).
Hence, for every \(\tau \in \Gal(K/\QQ)\) and every \(u \in K\), we have
\[
\begin{aligned}
\Norm_{K/F}\bigl(\tau(u)\bigr)
& = \prod_{\gamma \in \Gal(K/F)}\gamma\bigl(\tau(u)\bigr)\\
& = \tau\left(\prod_{\gamma \in \Gal(K/F)}(\tau^{-1}\gamma\tau)(u)\right)\\
& = \tau\bigl(\Norm_{K/F}(u)\bigr).
\end{aligned}
\]
That is,
\[
\Norm_{K/F} \circ \tau = \bigl(\tau|_{F}\bigr) \circ \Norm_{K/F}.
\]
For every \(\tau \in \Gal(K/\QQ)\) and every \(u \in \mathcal{E}(K/F)\), we have \(\tau(u) \in \mathcal{O}_{K}^{\times}\) and
\[
\Norm_{K/F}(\tau(u)) = \tau(\Norm_{K/F}(u)) = \tau(1) = 1.
\]
Thus \(\tau(u) \in \mathcal{E}(K/F)\).

In particular, for each \(a \in \{1, \dots, r(K) + 1\}\) and each \(j \in \{1, \dots, \rho\}\), we have
\[
\sigma_{a}(v_{j}),\ \kappa(\sigma_{a}(v_{j})) \in \mathcal{E}(K/F).
\]
Moreover,
\[
\sigma_{a}(v_{j})\kappa(\sigma_{a}(v_{j})) = |\sigma_{a}(v_{j})|^{2} > 0.
\]
These facts, together with the fact that \(\mathcal{E}(K/F)\) is closed under multiplication, give
\[
\bigl(\sigma_{a}(v_{j})\kappa(\sigma_{a}(v_{j}))\bigr)^{w_{a}} \in \mathcal{E}(K/F) \cap U(K) = U^{1}(K/F).
\]
Consequently, there exist integers \(c_{a,j,k}\) \((1 \le k \le e)\) such that
\[
2w_{a}\log|\sigma_{a}(v_{j})| = \sum_{k = 1}^{e}c_{a,j,k}y_{k}.
\]
\end{proof}

\begin{remark}\label{rem:relative-log-decomposition-factor-two}
We add a comment on the factor \(2\) that appears on the left-hand sides of the two expressions in Lemma~\ref{lem:relative-log-decomposition}. In what follows, we use the notation of Lemma~\ref{lem:relative-log-decomposition}.

As for the first expression, we can also express the number \(w_{a}\log|\sigma_{a}(\varepsilon_{j})|\), obtained from the left-hand side by removing the factor \(2\), as a linear combination of \(x_{1}, \dots, x_{\delta_{+}}\) with integer coefficients.
Indeed, since \(F/\QQ\) is a Galois extension, both \(\sigma_{a}(\varepsilon_{j})\) and its complex conjugate belong to \(\mathcal{O}_{F}^{\times}\).
If \(K\) is totally real, then \(w_{a} = 1\) and
\[
|\sigma_{a}(\varepsilon_{j})| = \pm\sigma_{a}(\varepsilon_{j}) \in U(F).
\]
If \(K\) is totally imaginary, then \(w_{a} = 2\) and
\[
|\sigma_{a}(\varepsilon_{j})|^{2} = \sigma_{a}(\varepsilon_{j})\overline{\sigma_{a}(\varepsilon_{j})} \in U(F).
\]

In contrast, as for the second expression, it is not always possible to express the number \(w_{a}\log|\sigma_{a}(v_{j})|\), obtained from the left-hand side by removing the factor \(2\), as a linear combination of \(y_{1}, \dots, y_{e}\) with integer coefficients.
We give such an example below.

Consider the cubic polynomial
\[
g(T) \coloneqq T^{3} + T^{2} - 2T - 1
\]
in an indeterminate \(T\).
Since \(g\) is a monic polynomial with integer coefficients and constant term \(-1\), its only possible rational roots are \(\pm 1\).
Since \(g(1) = -1\) and \(g(-1) = 1\), the polynomial \(g\) has no rational root, and since \(g\) is cubic, it is irreducible over \(\QQ\).
Since \(g(-1) = 1 > 0\) and \(g(0) = -1 < 0\), the intermediate value theorem shows that \(g\) has a negative real root \(\beta\).
The polynomial identity
\[
g(T^{2} - 2) = g(T)(T^{3} - T^{2} - 2T + 1)
\]
shows that \(\beta^{2} - 2\) is also a root of \(g\).
If \(\beta^{2} - 2 = \beta\), then \(\beta \in \{-1, 2\}\), which contradicts the fact that \(g\) has no rational root; hence \(\beta^{2} - 2 \ne \beta\).
By the relations between the roots and the coefficients of a polynomial, the sum of the three roots of \(g\) is \(-1\), and therefore the remaining root is
\[
-1 - \beta - (\beta^{2} - 2) = -\beta^{2} - \beta + 1.
\]

Put \(K \coloneqq \QQ(\beta)\) and \(F \coloneqq \QQ\).
Since all roots of \(g\) belong to \(K\), the field \(K\) is a totally real cubic field that is Galois over \(\QQ\).
Since \(r(K) = 2 > 0\), the extension \(K/F\) satisfies the assumptions of Lemma~\ref{lem:relative-log-decomposition}.
Moreover, \(\rho = r(K) - r(F) = 2\).
Since \(K \subseteq \RR\), we have \(\mu(K) = \{\pm 1\}\), and since \(\Norm_{K/F}(-1) = (-1)^{[K:F]} = -1\), we have \(\mathcal{E}(K/F) \cap \mu(K) = \{1\}\).
Hence \(v_{1}, v_{2}\) form a \(\ZZ\)-basis of \(\mathcal{E}(K/F)\).

Since \(g\) is a monic polynomial with integer coefficients, its three roots are algebraic integers.
Furthermore, since \(g\) has odd degree and constant term \(-1\), the product of its three roots is \(1\).
Therefore, \(\beta^{-1}\) equals the product of the remaining two roots and belongs to \(\mathcal{O}_{K}\).
Thus \(\beta \in \mathcal{O}_{K}^{\times}\).
In addition, since \(g\) is irreducible over \(\QQ\), its three roots are distinct and coincide with the images of \(\beta\) under the three elements of \(\Gal(K/\QQ)\).
Hence, by the product expression for the norm in a Galois extension stated in Section~\ref{sec:preliminaries}, \(\Norm_{K/F}(\beta)\) equals the product of the three roots of \(g\), and \(\Norm_{K/F}(\beta) = 1\).
Consequently, \(\beta \in \mathcal{E}(K/F)\).
Then there exist integers \(B_{1}, B_{2}\) such that
\[
\beta = v_{1}^{B_{1}}v_{2}^{B_{2}}.
\]
If both \(v_{1}\) and \(v_{2}\) are positive, then the right-hand side is positive, which contradicts \(\beta < 0\).
Hence there exists \(j_{0} \in \{1, 2\}\) such that \(v_{j_{0}} < 0\).

Since \(K\) is a totally real cubic field, the system of representatives \((\sigma_{1}, \sigma_{2}, \sigma_{3})\) of the embeddings consists of all embeddings of \(K\), and \(w_{a} = 1\) for each \(a \in \{1, 2, 3\}\).
The inclusion map \(K \hookrightarrow \CC\) coincides with one of \(\sigma_{1}, \sigma_{2}, \sigma_{3}\), and hence, after renumbering, we may assume without loss of generality that \(\sigma_{1}\) is the inclusion map \(K \hookrightarrow \CC\).
Then
\[
|\sigma_{1}(v_{j_{0}})|^{w_{1}} = |v_{j_{0}}| = -v_{j_{0}}.
\]
Since \(v_{j_{0}} \in \mathcal{E}(K/F)\) and \([K:F] = 3\), we find
\[
\Norm_{K/F}(-v_{j_{0}}) = (-1)^{[K:F]}\Norm_{K/F}(v_{j_{0}}) = -1.
\]
Thus
\[
|\sigma_{1}(v_{j_{0}})|^{w_{1}} \notin U^{1}(K/F).
\]

If there exist integers \(c_{1}, \dots, c_{e}\) such that
\[
w_{1}\log|\sigma_{1}(v_{j_{0}})| = \sum_{k = 1}^{e}c_{k}y_{k},
\]
then the application of the exponential function to both sides gives
\[
|\sigma_{1}(v_{j_{0}})|^{w_{1}} = \prod_{k = 1}^{e}\xi_{k}^{c_{k}} \in U^{1}(K/F),
\]
which is a contradiction.
Therefore, in this example, we cannot express the number \(w_{1}\log|\sigma_{1}(v_{j_{0}})|\) as a linear combination of \(y_{1}, \dots, y_{e}\) with integer coefficients.
\end{remark}

\begin{lemma}\label{lem:relative-regulator-polynomial}
Let \(F\) and \(K\) be number fields that are Galois over \(\QQ\), with \(F \subseteq K\) and \(r(K) > 0\).
Put
\[
\delta \coloneqq r(F),
\qquad \rho \coloneqq r(K) - r(F),
\qquad \delta_{+} \coloneqq r(F^{+}),
\]
and let \(e\) be the rank of \(U^{1}(K/F)\).
Let \(\eta_{1}, \dots, \eta_{\delta_{+}}\) be a \(\ZZ\)-basis of \(U(F)\),
let \(\xi_{1}, \dots, \xi_{e}\) be a \(\ZZ\)-basis of \(U^{1}(K/F)\), and put
\[
\begin{aligned}
x_{k}
& \coloneqq \log\eta_{k} \quad (1 \le k \le \delta_{+}),\\
y_{k}
& \coloneqq \log\xi_{k} \quad (1 \le k \le e).
\end{aligned}
\]
Moreover, let \(X_{1}, \dots, X_{\delta_{+}}, Y_{1}, \dots, Y_{e}\) be distinct indeterminates, and put
\[
\begin{aligned}
\boldsymbol{X}
& \coloneqq (X_{1}, \dots, X_{\delta_{+}}),
& \qquad \boldsymbol{x}
& \coloneqq (x_{1}, \dots, x_{\delta_{+}}),\\
\boldsymbol{Y}
& \coloneqq (Y_{1}, \dots, Y_{e}),
& \qquad \boldsymbol{y}
& \coloneqq (y_{1}, \dots, y_{e}).
\end{aligned}
\]
If \(\delta_{+} = 0\) or \(e = 0\), then we take the corresponding tuples to be empty tuples.
Then, with respect to \((\boldsymbol{X}, \boldsymbol{Y})\), there exists a nonzero multihomogeneous polynomial of multidegree
\[
(\delta,\rho),
\]
say
\[
P_{K/F} \in \QQ[\boldsymbol{X}, \boldsymbol{Y}],
\]
such that
\[
\reg(K) = \left|P_{K/F}(\boldsymbol{x}, \boldsymbol{y})\right|.
\]
\end{lemma}

\begin{proof}
Choose \(\varepsilon_{1}, \dots, \varepsilon_{\delta}\) and \(v_{1}, \dots, v_{\rho}\) as in Lemma~\ref{lem:relative-unit-kernel}.
Let
\[
I \coloneqq \left[\mathcal{O}_{K}^{\times}/\mu(K):\left\langle \varepsilon_{1}\mu(K), \dots, \varepsilon_{\delta}\mu(K), v_{1}\mu(K), \dots, v_{\rho}\mu(K) \right\rangle\right]
\]
be the index of the subgroup generated by the cosets of these units.

Let \((\sigma_{1}, \dots, \sigma_{r(K) + 1})\) be the system of representatives of the embeddings of \(K\) in the definition of the logarithmic matrix in Section~\ref{sec:preliminaries},
and let \(w_{1}, \dots, w_{r(K) + 1}\) be the corresponding weights.
By Lemma~\ref{lem:relative-log-decomposition}, there exist integers \(b_{a,j,k}\) and \(c_{a,j,k}\) such that
\[
2w_{a}\log|\sigma_{a}(\varepsilon_{j})| = \sum_{k = 1}^{\delta_{+}}b_{a,j,k}x_{k} \quad (1 \le a \le r(K) + 1,\ 1 \le j \le \delta)
\]
and
\[
2w_{a}\log|\sigma_{a}(v_{j})| = \sum_{k = 1}^{e}c_{a,j,k}y_{k} \quad (1 \le a \le r(K) + 1,\ 1 \le j \le \rho).
\]

We now put
\[
(u_{1}, \dots, u_{r(K)}) \coloneqq (\varepsilon_{1}, \dots, \varepsilon_{\delta}, v_{1}, \dots, v_{\rho}).
\]
In addition, for each \(a, j \in \{1, \dots, r(K)\}\), we put
\[
\mathcal{L}_{a,j}(\boldsymbol{X}, \boldsymbol{Y}) \coloneqq
\begin{cases}
\displaystyle \sum_{k = 1}^{\delta_{+}}b_{a,j,k}X_{k}, & 1 \le j \le \delta,\\
\displaystyle \sum_{k = 1}^{e}c_{a,j - \delta,k}Y_{k}, & \delta < j \le r(K).
\end{cases}
\]
Then
\[
2w_{a}\log|\sigma_{a}(u_{j})| = \mathcal{L}_{a,j}(\boldsymbol{x}, \boldsymbol{y}) \quad (1 \le a,j \le r(K)).
\]
Accordingly, we put
\[
P_{0}(\boldsymbol{X}, \boldsymbol{Y}) \coloneqq \det\bigl(\mathcal{L}_{a,j}(\boldsymbol{X}, \boldsymbol{Y})\bigr)_{1 \le a,j \le r(K)} \in \ZZ[\boldsymbol{X}, \boldsymbol{Y}].
\]
The first \(\delta\) columns of the matrix used in the definition of \(P_{0}\) consist only of homogeneous linear forms in \(\boldsymbol{X}\),
and the remaining \(\rho\) columns consist only of homogeneous linear forms in \(\boldsymbol{Y}\).
Therefore, each monomial of \(P_{0}\) has multidegree
\[
(\delta,\rho).
\]

The above equalities and the definition of \(P_{0}\) give
\[
P_{0}(\boldsymbol{x}, \boldsymbol{y}) = 2^{r(K)}\det\mathcal{M}_{K}^{(r(K) + 1)}(u_{1}, \dots, u_{r(K)}).
\]
By the index formula~\eqref{eq:index-formula} for logarithmic matrices, we have
\[
\left|\det\mathcal{M}_{K}^{(r(K) + 1)}(u_{1}, \dots, u_{r(K)})\right| = I\reg(K).
\]
Consequently,
\[
\left|P_{0}(\boldsymbol{x}, \boldsymbol{y})\right| = 2^{r(K)}I\reg(K).
\]
It therefore suffices to put
\[
P_{K/F} \coloneqq \frac{1}{2^{r(K)}I}P_{0}.
\]
Since \(\reg(K) > 0\), the polynomial \(P_{K/F}\) is nonzero.
Moreover, as shown above, each monomial of \(P_{0}\) has multidegree \((\delta,\rho)\),
and hence \(P_{K/F}\) is a multihomogeneous polynomial of multidegree \((\delta,\rho)\) with respect to \((\boldsymbol{X}, \boldsymbol{Y})\).
\end{proof}

\section{Linear independence of logarithms}\label{sec:log-linear-independence}

In this section, we show that the logarithms substituted in the polynomial expressions of Section~\ref{sec:regulator-polynomial} are linearly independent over \(\QQ\).

We first treat finite families of number fields in which the unit group of the intersection of any two distinct fields has rank \(0\).
For this purpose, we prepare the following lemma.

\begin{lemma}\label{lem:norm-through-intersection}
Let \(K\) and \(K'\) be number fields that are Galois over \(\QQ\), and put
\[
F \coloneqq K \cap K'.
\]
Moreover, let \(E\) be a number field that satisfies \(KK' \subseteq E\).
Then, for every \(u \in K\),
\[
\Norm_{E/K'}(u) = \Norm_{K/F}(u)^{[E:KK']}.
\]
Furthermore, if
\[
u \in \mathcal{O}_{K}^{\times} \cap \RR,
\]
then
\[
\Norm_{K/F}(u) \in \mathcal{O}_{F}^{\times} \cap \RR.
\]
\end{lemma}

\begin{proof}
Since \(K/\QQ\) is a Galois extension, \(KK'/K'\) and \(K/F\) are also Galois extensions, and the restriction map
\[
\Gal(KK'/K') \longrightarrow \Gal(K/F),
\qquad \sigma \longmapsto \sigma|_{K}
\]
is an isomorphism~\cite[Proposition~3.19]{MilneFT}.
Hence, for \(u \in K\), the product expression for the norm in a Galois extension stated in Section~\ref{sec:preliminaries} gives
\[
\begin{aligned}
\Norm_{KK'/K'}(u)
& = \prod_{\sigma \in \Gal(KK'/K')} \sigma(u)\\
& = \prod_{\sigma \in \Gal(KK'/K')} \bigl(\sigma|_{K}\bigr)(u)\\
& = \prod_{\gamma \in \Gal(K/F)} \gamma(u)\\
& = \Norm_{K/F}(u).
\end{aligned}
\]
Here, the third equality uses the fact that, since the above restriction map is an isomorphism,
\(\sigma|_{K}\) runs through all elements of \(\Gal(K/F)\) exactly once as \(\sigma\) runs through all elements of \(\Gal(KK'/K')\).
Since \(u \in KK'\), we have
\[
\Norm_{E/KK'}(u) = u^{[E:KK']}.
\]
Therefore, the transitivity and the multiplicativity of the norm yield
\[
\begin{aligned}
\Norm_{E/K'}(u)
& = \Norm_{KK'/K'}\bigl(\Norm_{E/KK'}(u)\bigr)\\
& = \Norm_{KK'/K'}\bigl(u^{[E:KK']}\bigr)\\
& = \Norm_{KK'/K'}(u)^{[E:KK']}\\
& = \Norm_{K/F}(u)^{[E:KK']}.
\end{aligned}
\]

Next, suppose that \(u \in \mathcal{O}_{K}^{\times} \cap \RR\).
As stated in Section~\ref{sec:preliminaries}, \(\Norm_{K/F}(\mathcal{O}_{K}^{\times}) \subseteq \mathcal{O}_{F}^{\times}\), and hence
\[
\Norm_{K/F}(u) \in \mathcal{O}_{F}^{\times}.
\]
Since \(K/\QQ\) is a Galois extension, the restriction of complex conjugation
\[
\kappa \colon K \longrightarrow K,
\qquad x \longmapsto \overline{x}
\]
is well-defined, and \(\kappa \in \Gal(K/\QQ)\).
Since both \(K/\QQ\) and \(K'/\QQ\) are Galois extensions and \(F = K \cap K'\), the extension \(F/\QQ\) is also a Galois extension.
Therefore, \(\Gal(K/F)\) is a normal subgroup of \(\Gal(K/\QQ)\). Thus the map
\[
\Gal(K/F) \longrightarrow \Gal(K/F),
\qquad \gamma \longmapsto \kappa\gamma\kappa^{-1}
\]
is well-defined. This map is a bijection; indeed, \(\gamma \mapsto \kappa^{-1}\gamma\kappa\) gives its inverse. Consequently,
\[
\begin{aligned}
\kappa\bigl(\Norm_{K/F}(u)\bigr)
& = \prod_{\gamma \in \Gal(K/F)} \kappa\gamma(u)\\
& = \prod_{\gamma \in \Gal(K/F)} \bigl(\kappa\gamma\kappa^{-1}\bigr)\bigl(\kappa(u)\bigr)\\
& = \Norm_{K/F}\bigl(\kappa(u)\bigr).
\end{aligned}
\]
Since \(u \in \RR\) gives \(\kappa(u) = u\), we obtain
\[
\kappa\bigl(\Norm_{K/F}(u)\bigr) = \Norm_{K/F}(u).
\]
It follows that
\[
\Norm_{K/F}(u) \in \RR.
\]
We conclude that
\[
\Norm_{K/F}(u) \in \mathcal{O}_{F}^{\times} \cap \RR.
\]
\end{proof}

\begin{lemma}\label{lem:linear-independence}
Let \(m \ge 1\) be an integer, and let \(K_{1}, \dots, K_{m}\) be number fields that are Galois over \(\QQ\).
Assume that
\[
r(K_{i} \cap K_{j}) = 0
\]
for all distinct \(i, j \in \{1, \dots, m\}\).
For each \(i \in \{1, \dots, m\}\), put \(r_{i,+} \coloneqq r(K_{i}^{+})\),
and let \(\eta_{i, 1}, \dots, \eta_{i, r_{i,+}}\) be a \(\ZZ\)-basis of \(U(K_{i})\) (we take the empty tuple when \(r_{i,+} = 0\)).
Furthermore, put
\[
x_{i, k} \coloneqq \log\eta_{i, k} \quad (1 \le i \le m,\ 1 \le k \le r_{i,+}).
\]
Then, if \(r_{1,+} + \dots + r_{m,+} \ge 1\), these \(r_{1,+} + \dots + r_{m,+}\) real numbers \(x_{i, k}\) are linearly independent over \(\QQ\).
\end{lemma}

\begin{proof}
By clearing denominators, we can reduce a linear relation with rational coefficients to a linear relation with integer coefficients.
Accordingly, take arbitrary integers \(A_{i, k}\) indexed by \(1 \le i \le m\) and \(1 \le k \le r_{i,+}\), and assume that
\[
\sum_{i = 1}^{m}\sum_{k = 1}^{r_{i,+}}A_{i, k}x_{i, k} = 0.
\]
It suffices to show that \(A_{i, k} = 0\) for all \(i \in \{1, \dots, m\}\) and all \(k \in \{1, \dots, r_{i,+}\}\).
For each \(i \in \{1, \dots, m\}\), put
\[
u_{i} \coloneqq \prod_{k = 1}^{r_{i,+}}\eta_{i, k}^{A_{i, k}} \in U(K_{i}),
\]
where, when \(r_{i,+} = 0\), we put \(u_{i} \coloneqq 1\) as the empty product.
If we apply the exponential function to both sides of the assumed linear relation, then we obtain
\[
\begin{aligned}
1
& = \exp\left(\sum_{i = 1}^{m}\sum_{k = 1}^{r_{i,+}}A_{i, k}x_{i, k}\right)\\
& = \prod_{i = 1}^{m}\prod_{k = 1}^{r_{i,+}}\exp\bigl(A_{i, k}\log\eta_{i, k}\bigr)\\
& = \prod_{i = 1}^{m}\prod_{k = 1}^{r_{i,+}}\eta_{i, k}^{A_{i, k}}\\
& = \prod_{i = 1}^{m}u_{i}.
\end{aligned}
\]

Put \(\widetilde{K} \coloneqq K_{1}\cdots K_{m}\).
Fix an arbitrary \(j \in \{1, \dots, m\}\).
Since \(u_{j} \in K_{j}\), we have \(\Norm_{\widetilde{K}/K_{j}}(u_{j}) = u_{j}^{[\widetilde{K}:K_{j}]}\).
Hence, if we apply the norm \(\Norm_{\widetilde{K}/K_{j}}\) to the equality \(\prod_{i = 1}^{m}u_{i} = 1\) and use the multiplicativity of the norm, then we get
\[
1 = u_{j}^{[\widetilde{K}:K_{j}]} \prod_{\substack{1 \le i \le m \\ i \ne j}} \Norm_{\widetilde{K}/K_{j}}(u_{i}).
\]
Here, when \(m = 1\), we take the product over \(i\) on the right-hand side to be \(1\) as the empty product.

Take an arbitrary \(i \in \{1, \dots, m\} \setminus \{j\}\), and put
\[
F_{i, j} \coloneqq K_{i} \cap K_{j}.
\]
Since \(K_{i}K_{j} \subseteq \widetilde{K}\) and
\[
u_{i} \in U(K_{i}) \subseteq \mathcal{O}_{K_{i}}^{\times} \cap \RR,
\]
Lemma~\ref{lem:norm-through-intersection} applied with
\[
K = K_{i},
\qquad K' = K_{j},
\qquad E = \widetilde{K}
\]
shows that
\[
\Norm_{K_{i}/F_{i, j}}(u_{i}) \in \mathcal{O}_{F_{i, j}}^{\times} \cap \RR
\]
and
\[
\Norm_{\widetilde{K}/K_{j}}(u_{i}) = \Norm_{K_{i}/F_{i, j}}(u_{i})^{[\widetilde{K}:K_{i}K_{j}]}.
\]

On the other hand, since \(r(F_{i, j}) = 0\) by assumption, \(\mathcal{O}_{F_{i, j}}^{\times}/\mu(F_{i, j})\) is a free abelian group of rank \(0\), that is, the trivial group, and
\[
\mathcal{O}_{F_{i, j}}^{\times} = \mu(F_{i, j}).
\]
Hence
\[
\Norm_{K_{i}/F_{i, j}}(u_{i}) \in \mu(F_{i, j}) \cap \RR = \{\pm 1\}.
\]
Since \(i\) was arbitrary, we obtain
\[
\Norm_{\widetilde{K}/K_{j}}(u_{i}) \in \{\pm 1\}
\]
for all \(i \in \{1, \dots, m\} \setminus \{j\}\).
If we use this in the above product equality, then we find
\[
u_{j}^{[\widetilde{K}:K_{j}]} \in \{\pm 1\}.
\]
On the other hand, since \(u_{j} \in U(K_{j}) \subseteq \RR_{ > 0}\), the left-hand side is positive, and thus
\[
u_{j}^{[\widetilde{K}:K_{j}]} = 1.
\]
Therefore, \(u_{j} = 1\).
Since \(j\) was arbitrary, we have \(u_{j} = 1\) for all \(j \in \{1, \dots, m\}\).
For each \(i \in \{1, \dots, m\}\), \(\eta_{i, 1}, \dots, \eta_{i, r_{i,+}}\) is a \(\ZZ\)-basis of \(U(K_{i})\),
and consequently \(A_{i, k} = 0\) for all \(i \in \{1, \dots, m\}\) and all \(k \in \{1, \dots, r_{i,+}\}\).
\end{proof}

\begin{remark}\label{rem:linear-independence-small-m}
In the cases \(m = 1\) and \(m = 2\) of Lemma~\ref{lem:linear-independence}, there is an alternative proof that does not use the Galois property of \(K_{i}/\QQ\).
In what follows, we use the notation of the proof of Lemma~\ref{lem:linear-independence}.
When \(m = 1\), the equality \(\prod_{i = 1}^{m}u_{i} = 1\) gives \(u_{1} = 1\).
When \(m = 2\), put \(F \coloneqq K_{1} \cap K_{2}\).
The equality \(u_{1}u_{2} = 1\) gives
\[
u_{1} = u_{2}^{-1} \in F \cap \RR_{ > 0}.
\]
In addition, since both \(u_{1}\) and \(u_{1}^{-1}\) are algebraic integers, we see that \(u_{1} \in U(F)\).
As stated in Section~\ref{sec:maximal-real-subfield}, \(U(F)\) is a free abelian group of rank \(r(F^{+})\).
Since \(F^{+} \subseteq F\), Lemma~\ref{lem:relative-rank} and the assumption \(r(F) = 0\) give \(r(F^{+}) = 0\).
Therefore, \(U(F)\) is a free abelian group of rank \(0\), that is, the trivial group, and \(u_{1} = 1\).
Hence \(u_{2} = 1\) as well.
In either case, we have \(u_{i} = 1\) for all \(i \in \{1, \dots, m\}\).
Since, for each \(i \in \{1, \dots, m\}\), \(\eta_{i, 1}, \dots, \eta_{i, r_{i,+}}\) is a \(\ZZ\)-basis of \(U(K_{i})\), all \(A_{i, k}\) are \(0\).

In contrast, when \(m \ge 3\), the equality
\[
u_{1}\cdots u_{m} = 1
\]
implies, for example,
\[
u_{1} = (u_{2}\cdots u_{m})^{-1} \in K_{1} \cap (K_{2}\cdots K_{m}),
\]
but the condition \(r(K_{i} \cap K_{j}) = 0\) for all distinct \(i, j\) alone does not imply
\[
r\bigl(K_{1} \cap (K_{2}\cdots K_{m})\bigr) = 0
\]
in general. Indeed, put
\[
K_{1} \coloneqq \QQ(\sqrt{2}),
\qquad K_{2} \coloneqq \QQ(\sqrt{3}),
\qquad K_{3} \coloneqq \QQ(\sqrt{6}).
\]
These are distinct quadratic fields, and hence
\[
K_{i} \cap K_{j} = \QQ
\]
for all distinct \(i, j\). On the other hand, \(\sqrt{2} = \sqrt{6}/\sqrt{3}\) gives
\[
K_{1} \subseteq K_{2}K_{3},
\]
and hence
\[
K_{1} \cap (K_{2}K_{3}) = K_{1};
\]
thus the unit group of this intersection has rank \(1\).
Therefore, the above argument for the case \(m = 2\) does not extend directly to the case \(m \ge 3\).
In the proof of Lemma~\ref{lem:linear-independence}, we apply Lemma~\ref{lem:norm-through-intersection} for each \(i \ne j\)
and express \(\Norm_{\widetilde{K}/K_{j}}(u_{i})\) as a power of \(\Norm_{K_{i}/F_{i, j}}(u_{i})\).
Here, \(F_{i, j} = K_{i} \cap K_{j}\) is the intersection of two fields, and the use of the assumption \(r(F_{i, j}) = 0\) on this field avoids the problem described above.
Lemma~\ref{lem:norm-through-intersection} uses the Galois property of \(K_{i}/\QQ\) and of \(K_{j}/\QQ\).
\end{remark}

Next, for two number fields \(K_{1}, K_{2}\) and their intersection \(F\), we treat the logarithms of \(\ZZ\)-bases of \(U(F)\), \(U^{1}(K_{1}/F)\), and \(U^{1}(K_{2}/F)\).
The following lemma holds for general number fields, but in this paper we use it only in the case where both \(K_{1}\) and \(K_{2}\) are Galois over \(\QQ\).

\begin{lemma}\label{lem:relative-log-independence}
Let \(K_{1}\) and \(K_{2}\) be number fields, and put
\[
F \coloneqq K_{1} \cap K_{2}.
\]
Put \(\delta_{+} \coloneqq r(F^{+})\),
and let \(e_{i}\) be the rank of \(U^{1}(K_{i}/F)\) \((i = 1, 2)\).
Let \(\eta_{1}, \dots, \eta_{\delta_{+}}\) be a \(\ZZ\)-basis of \(U(F)\),
and let \(\xi_{i,1}, \dots, \xi_{i,e_{i}}\) be a \(\ZZ\)-basis of \(U^{1}(K_{i}/F)\) \((i = 1, 2)\).
For each of these three tuples, if the corresponding rank is \(0\), then we take the tuple to be the empty tuple.
Then, if \(\delta_{+} + e_{1} + e_{2} \ge 1\), the real numbers
\[
\log\eta_{1}, \dots, \log\eta_{\delta_{+}}, \log\xi_{1,1}, \dots, \log\xi_{1,e_{1}}, \log\xi_{2,1}, \dots, \log\xi_{2,e_{2}}
\]
are linearly independent over \(\QQ\).
\end{lemma}

\begin{proof}
Assume that \(\delta_{+} + e_{1} + e_{2} \ge 1\).
By clearing denominators, we can reduce a linear relation with rational coefficients to a linear relation with integer coefficients.
Accordingly, take arbitrary
\(A_{k} \in \ZZ\) \((1 \le k \le \delta_{+})\),
\(B_{k} \in \ZZ\) \((1 \le k \le e_{1})\), and
\(C_{k} \in \ZZ\) \((1 \le k \le e_{2})\),
and assume that
\[
\sum_{k = 1}^{\delta_{+}}A_{k}\log\eta_{k} + \sum_{k = 1}^{e_{1}}B_{k}\log\xi_{1,k} + \sum_{k = 1}^{e_{2}}C_{k}\log\xi_{2,k} = 0.
\]
Put
\[
a \coloneqq \prod_{k = 1}^{\delta_{+}}\eta_{k}^{A_{k}},
\qquad b \coloneqq \prod_{k = 1}^{e_{1}}\xi_{1,k}^{B_{k}},
\qquad c \coloneqq \prod_{k = 1}^{e_{2}}\xi_{2,k}^{C_{k}},
\]
where we take the product corresponding to an empty tuple to be \(1\).
If we apply the exponential function to both sides of the assumed linear relation, then we obtain
\[
abc = 1.
\]
Moreover,
\[
a \in U(F),
\qquad b \in U^{1}(K_{1}/F),
\qquad c \in U^{1}(K_{2}/F).
\]

The equality \(c^{-1} = ab\) gives \(c \in K_{1}\), and the definition of \(c\) gives \(c \in K_{2}\) as well; hence
\[
c \in K_{1} \cap K_{2} = F.
\]
It follows that
\[
1 = \Norm_{K_{2}/F}(c) = c^{[K_{2}:F]}.
\]
Since \(c > 0\), we conclude that \(c = 1\).
Hence \(ab = 1\).
Since \(a \in U(F) \subseteq F\), we have \(b = a^{-1} \in F\). Thus
\[
1 = \Norm_{K_{1}/F}(b) = b^{[K_{1}:F]}.
\]
Since \(b > 0\), we obtain \(b = 1\), and therefore \(a = 1\).
Since each tuple is a \(\ZZ\)-basis of a free abelian group, all \(A_{k}\), \(B_{k}\), and \(C_{k}\) are \(0\).
\end{proof}

\section{Algebraic independence of polynomials}\label{sec:polynomial-independence}

In this section, we prove the algebraic independence of homogeneous polynomials and of multihomogeneous polynomials by a comparison of multidegrees.

We first treat homogeneous polynomials that have no indeterminates in common.

\begin{lemma}\label{lem:disjoint-homogeneous}
Let \(L\) be a field, and let \(m \ge 1\) be an integer.
For each \(i \in \{1, \dots, m\}\), let \(d_{i} \ge 1\) be an integer, and let
\[
\boldsymbol{X}_{i} = (X_{i, 1}, \dots, X_{i, d_{i}})
\]
be a tuple of indeterminates.
We assume that, for all distinct \(i, j \in \{1, \dots, m\}\), the tuples \(\boldsymbol{X}_{i}\) and \(\boldsymbol{X}_{j}\) have no component in common.
For each \(i \in \{1, \dots, m\}\), let \(r_{i} \ge 1\) be an integer, and let
\[
P_{i} \in L[\boldsymbol{X}_{i}] = L[X_{i, 1}, \dots, X_{i, d_{i}}]
\]
be a nonzero homogeneous polynomial of total degree \(r_{i}\).
Then the polynomials
\[
P_{1}(\boldsymbol{X}_{1}), \dots, P_{m}(\boldsymbol{X}_{m})
\]
are algebraically independent over \(L\).
\end{lemma}

\begin{proof}
Take an arbitrary nonzero polynomial
\[
H(T_{1}, \dots, T_{m}) = \sum_{\boldsymbol{\nu} \in S} \lambda_{\boldsymbol{\nu}} \prod_{i = 1}^{m}T_{i}^{\nu_{i}} \in L[T_{1}, \dots, T_{m}]
\]
in indeterminates \(T_{1}, \dots, T_{m}\).
Here, \(S\) is a nonempty finite subset of \(\ZZ_{ \ge 0}^{m}\),
and \(\lambda_{\boldsymbol{\nu}} \in L \setminus \{0\}\) for each \(\boldsymbol{\nu} = (\nu_{1}, \dots, \nu_{m}) \in S\).
Then
\[
H\bigl( P_{1}(\boldsymbol{X}_{1}), \dots, P_{m}(\boldsymbol{X}_{m}) \bigr) = \sum_{\boldsymbol{\nu} \in S} \lambda_{\boldsymbol{\nu}} \prod_{i = 1}^{m} P_{i}(\boldsymbol{X}_{i})^{\nu_{i}}.
\]

For each \(\boldsymbol{\nu} \in S\), the polynomial
\[
\lambda_{\boldsymbol{\nu}} \prod_{i = 1}^{m} P_{i}(\boldsymbol{X}_{i})^{\nu_{i}}
\]
is nonzero because the polynomial ring \(L[\boldsymbol{X}_{1}, \dots, \boldsymbol{X}_{m}]\) over the field \(L\) is an integral domain,
and the multidegree with respect to \((\boldsymbol{X}_{1}, \dots, \boldsymbol{X}_{m})\) of each monomial appearing in this polynomial equals
\[
(\nu_{1}r_{1}, \dots, \nu_{m}r_{m}).
\]
Since all \(r_{i}\) are positive, distinct \(\boldsymbol{\nu}\) give distinct multidegrees.
Hence the terms corresponding to distinct \(\boldsymbol{\nu}\) do not cancel one another.
Consequently,
\[
H\bigl( P_{1}(\boldsymbol{X}_{1}), \dots, P_{m}(\boldsymbol{X}_{m}) \bigr) \ne 0.
\]
Since \(H\) was an arbitrary nonzero polynomial,
\(P_{1}(\boldsymbol{X}_{1}), \dots, P_{m}(\boldsymbol{X}_{m})\)
are algebraically independent over \(L\).
\end{proof}

Next, we treat two multihomogeneous polynomials that share a tuple of indeterminates.

\begin{lemma}\label{lem:shared-homogeneous-pair}
Let \(L\) be a field, and let \(d_{X}, d_{Y}, d_{Z}\) be nonnegative integers.
Let \(X_{1}, \dots, X_{d_{X}}, Y_{1}, \dots, Y_{d_{Y}}, Z_{1}, \dots, Z_{d_{Z}}\) be distinct indeterminates, and put
\[
\begin{aligned}
\boldsymbol{X}
& \coloneqq (X_{1}, \dots, X_{d_{X}}),\\
\boldsymbol{Y}
& \coloneqq (Y_{1}, \dots, Y_{d_{Y}}),\\
\boldsymbol{Z}
& \coloneqq (Z_{1}, \dots, Z_{d_{Z}}).
\end{aligned}
\]
If \(d_{X} = 0\), \(d_{Y} = 0\), or \(d_{Z} = 0\), then we take the corresponding tuple to be the empty tuple.
Moreover, let \(\delta, \rho_{1}, \rho_{2}\) be nonnegative integers such that
\[
\delta + \rho_{1} > 0,
\qquad \delta + \rho_{2} > 0,
\qquad \rho_{1} + \rho_{2} > 0.
\]
Furthermore, let
\[
P_{1} \in L[\boldsymbol{X}, \boldsymbol{Y}],
\qquad P_{2} \in L[\boldsymbol{X}, \boldsymbol{Z}]
\]
be nonzero polynomials, and regard them as elements of \(L[\boldsymbol{X}, \boldsymbol{Y}, \boldsymbol{Z}]\).
Suppose that \(P_{1}\) is a multihomogeneous polynomial of multidegree
\[
(\delta, \rho_{1}, 0)
\]
with respect to \((\boldsymbol{X}, \boldsymbol{Y}, \boldsymbol{Z})\) and that \(P_{2}\) is a multihomogeneous polynomial of multidegree
\[
(\delta, 0, \rho_{2}).
\]
Then \(P_{1}\) and \(P_{2}\) are algebraically independent over \(L\).
\end{lemma}

\begin{proof}
Take an arbitrary nonzero polynomial
\[
H(T_{1}, T_{2}) = \sum_{(a,b) \in S}\lambda_{a,b}T_{1}^{a}T_{2}^{b} \in L[T_{1}, T_{2}]
\]
in indeterminates \(T_{1}, T_{2}\).
Here, \(S\) is a nonempty finite subset of \(\ZZ_{ \ge 0}^{2}\),
and \(\lambda_{a,b} \in L \setminus \{0\}\) for each \((a,b) \in S\).
The multidegree of each monomial appearing in each term
\[
\lambda_{a,b}P_{1}^{a}P_{2}^{b}
\]
of \(H(P_{1},P_{2})\) equals
\[
\bigl(\delta(a + b), \rho_{1}a, \rho_{2}b\bigr).
\]

The map
\[
\Phi \colon \ZZ_{ \ge 0}^{2} \longrightarrow \ZZ_{ \ge 0}^{3},
\qquad (a,b) \longmapsto \bigl(\delta(a + b), \rho_{1}a, \rho_{2}b\bigr)
\]
is injective.
Indeed, suppose that \((a,b), (a',b') \in \ZZ_{ \ge 0}^{2}\) satisfy \(\Phi(a,b) = \Phi(a',b')\);
that is, suppose that
\[
\bigl(\delta(a + b), \rho_{1}a, \rho_{2}b\bigr) = \bigl(\delta(a' + b'), \rho_{1}a', \rho_{2}b'\bigr).
\]

If \(\rho_{1} > 0\), then a comparison of the second components gives \(a = a'\).
In this case, if \(\rho_{2} > 0\), then a comparison of the third components gives \(b = b'\);
if \(\rho_{2} = 0\), then the assumption \(\delta + \rho_{2} > 0\) gives \(\delta > 0\),
and hence a comparison of the first components gives \(a + b = a' + b'\), which together with \(a = a'\) yields \(b = b'\).

If \(\rho_{1} = 0\), then the assumption \(\rho_{1} + \rho_{2} > 0\) gives \(\rho_{2} > 0\), and a comparison of the third components gives \(b = b'\).
Furthermore, the assumption \(\delta + \rho_{1} > 0\) gives \(\delta > 0\),
and thus a comparison of the first components gives \(a + b = a' + b'\), which together with \(b = b'\) yields \(a = a'\).

Consequently, \((a,b) = (a',b')\) in either case.
Therefore, \(\Phi\) is injective.

Each term \(\lambda_{a,b}P_{1}^{a}P_{2}^{b}\) is nonzero
because the polynomial ring \(L[\boldsymbol{X}, \boldsymbol{Y}, \boldsymbol{Z}]\) is an integral domain.
Moreover, by the injectivity of \(\Phi\), the terms corresponding to distinct \((a,b)\) do not cancel one another.
Hence \(H(P_{1},P_{2}) \ne 0\).
Since \(H\) was an arbitrary nonzero polynomial, we obtain the conclusion.
\end{proof}

\section{Proofs of the two main theorems and their corollaries}\label{sec:main}

In this section, we prove the two main theorems and deduce the two corollaries from them.

\begin{proof}[Proof of Theorem~\ref{thm:main-algebraic-independence}]
For each \(i \in \{1, \dots, m\}\), put
\[
r_{i} \coloneqq r(K_{i}),
\qquad r_{i,+} \coloneqq r(K_{i}^{+}).
\]
The assumption gives \(r_{i} > 0\), and Lemma~\ref{lem:rank-zero}~\textup{(iii)} gives \(r_{i,+} > 0\).
Let \(\eta_{i, 1}, \dots, \eta_{i, r_{i,+}}\) be a \(\ZZ\)-basis of \(U(K_{i})\), and put
\[
x_{i, k} \coloneqq \log\eta_{i, k} \quad (1 \le i \le m,\ 1 \le k \le r_{i,+}).
\]
In addition, put
\[
n \coloneqq \sum_{i = 1}^{m} r_{i,+}.
\]
Since all \(r_{i,+}\) are positive, we see that \(n \ge 1\).

By Lemma~\ref{lem:linear-independence}, the \(n\) real numbers \(x_{i, k}\) \((1 \le i \le m,\ 1 \le k \le r_{i,+})\) are linearly independent over \(\QQ\).
Since each \(\eta_{i, k}\) is an algebraic number that is a positive real number, WSC implies that these \(n\) numbers \(x_{i, k}\) are algebraically independent over \(\QQ\).

In what follows, for each \(i \in \{1, \dots, m\}\), we write
\[
\begin{aligned}
\boldsymbol{x}_{i}
& \coloneqq (x_{i, 1}, \dots, x_{i, r_{i,+}}),\\
\boldsymbol{X}_{i}
& \coloneqq (X_{i, 1}, \dots, X_{i, r_{i,+}})
\end{aligned}
\]
for a tuple of real numbers and a tuple of indeterminates, respectively.
We assume that, for all distinct \(i, j \in \{1, \dots, m\}\), the tuples \(\boldsymbol{X}_{i}\) and \(\boldsymbol{X}_{j}\) have no component in common.
For each \(i \in \{1, \dots, m\}\), we apply Lemma~\ref{lem:polynomial} with \(K = K_{i}\) and with the \(\ZZ\)-basis \(\eta_{i, 1}, \dots, \eta_{i, r_{i,+}}\).
In this application, we rename the indeterminates \(X_{1}, \dots, X_{r_{i,+}}\) of Lemma~\ref{lem:polynomial} as \(X_{i, 1}, \dots, X_{i, r_{i,+}}\).
Then there exists a nonzero homogeneous polynomial
\[
P_{i} \in \ZZ[X_{i, 1}, \dots, X_{i, r_{i,+}}]
\]
of total degree \(r_{i}\) with integer coefficients such that
\[
\reg(K_{i}) = |P_{i}(\boldsymbol{x}_{i})|.
\]
Since multiplication by \(-1\) preserves both the property of being nonzero and the total degree of a polynomial, we may assume, after changing the sign of \(P_{i}\) if necessary, that
\[
\reg(K_{i}) = P_{i}(\boldsymbol{x}_{i}).
\]

By Lemma~\ref{lem:disjoint-homogeneous}, the polynomials
\[
P_{1}(\boldsymbol{X}_{1}), \dots, P_{m}(\boldsymbol{X}_{m})
\]
are algebraically independent over \(\QQ\).
Hence, for every nonzero polynomial
\[
H \in \QQ[T_{1}, \dots, T_{m}]
\]
in indeterminates \(T_{1}, \dots, T_{m}\), the polynomial
\[
H\bigl( P_{1}(\boldsymbol{X}_{1}), \dots, P_{m}(\boldsymbol{X}_{m}) \bigr)
\]
is nonzero.

Since the components of all \(\boldsymbol{x}_{i}\), that is, the numbers \(x_{i, k}\) \((1 \le i \le m,\ 1 \le k \le r_{i,+})\), are algebraically independent over \(\QQ\), we have
\[
H\bigl(\reg(K_{1}), \dots, \reg(K_{m})\bigr) = H\bigl( P_{1}(\boldsymbol{x}_{1}), \dots, P_{m}(\boldsymbol{x}_{m}) \bigr) \ne 0.
\]
Since \(H\) was an arbitrary nonzero polynomial, we obtain the conclusion.
\end{proof}

\begin{proof}[Proof of Theorem~\ref{thm:galois-pair-criterion}]
Let \(T_{1}, T_{2}\) be indeterminates.

We first prove that \textup{(i)} implies \textup{(ii)}.
Assume \textup{(i)}, and put
\[
F \coloneqq K_{1} \cap K_{2}.
\]
Since \(K_{1}/\QQ\) and \(K_{2}/\QQ\) are Galois extensions, \(F/\QQ\) is also a Galois extension.
Put
\[
\delta \coloneqq r(F),
\qquad \rho_{i} \coloneqq r(K_{i}) - r(F) \quad (i = 1, 2).
\]
Lemma~\ref{lem:relative-rank} gives \(\rho_{1}, \rho_{2} \ge 0\).
Moreover, \textup{(i)} gives \(\delta + \rho_{i} = r(K_{i}) > 0\) \((i = 1, 2)\).
Furthermore, if \(\rho_{1} = \rho_{2} = 0\), then Lemma~\ref{lem:relative-rank} gives \(K_{1}^{+} = F^{+} = K_{2}^{+}\), which contradicts \textup{(i)}; hence \(\rho_{1} + \rho_{2} > 0\).
Consequently, we have
\begin{equation}\label{eq:multidegree-inequalities}
\delta + \rho_{1} > 0,
\qquad \delta + \rho_{2} > 0,
\qquad \rho_{1} + \rho_{2} > 0.
\end{equation}

Put \(\delta_{+} \coloneqq r(F^{+})\), and let \(e_{i}\) be the rank of \(U^{1}(K_{i}/F)\) \((i = 1, 2)\).
Choose a \(\ZZ\)-basis \(\eta_{1}, \dots, \eta_{\delta_{+}}\) of \(U(F)\),
a \(\ZZ\)-basis \(\xi_{1,1}, \dots, \xi_{1,e_{1}}\) of \(U^{1}(K_{1}/F)\),
and a \(\ZZ\)-basis \(\xi_{2,1}, \dots, \xi_{2,e_{2}}\) of \(U^{1}(K_{2}/F)\).
For each of these three tuples, if the corresponding rank is \(0\), then we take the tuple to be the empty tuple.
Moreover, put
\[
\begin{aligned}
x_{k}
& \coloneqq \log\eta_{k} \quad (1 \le k \le \delta_{+}),\\
y_{k}
& \coloneqq \log\xi_{1,k} \quad (1 \le k \le e_{1}),\\
z_{k}
& \coloneqq \log\xi_{2,k} \quad (1 \le k \le e_{2}),
\end{aligned}
\]
let \(X_{1}, \dots, X_{\delta_{+}}, Y_{1}, \dots, Y_{e_{1}}, Z_{1}, \dots, Z_{e_{2}}\) be distinct indeterminates, and put
\[
\begin{aligned}
\boldsymbol{X}
& \coloneqq (X_{1}, \dots, X_{\delta_{+}}),
& \qquad \boldsymbol{x}
& \coloneqq (x_{1}, \dots, x_{\delta_{+}}),\\
\boldsymbol{Y}
& \coloneqq (Y_{1}, \dots, Y_{e_{1}}),
& \qquad \boldsymbol{y}
& \coloneqq (y_{1}, \dots, y_{e_{1}}),\\
\boldsymbol{Z}
& \coloneqq (Z_{1}, \dots, Z_{e_{2}}),
& \qquad \boldsymbol{z}
& \coloneqq (z_{1}, \dots, z_{e_{2}}).
\end{aligned}
\]

For each \(i \in \{1, 2\}\), we apply Lemma~\ref{lem:relative-regulator-polynomial} with \(F = K_{1} \cap K_{2}\) and \(K = K_{i}\).
In this application, we read \(e\), \(\rho\), and \(\xi_{1}, \dots, \xi_{e}\) of Lemma~\ref{lem:relative-regulator-polynomial} as \(e_{i}\), \(\rho_{i}\), and \(\xi_{i,1}, \dots, \xi_{i,e_{i}}\), respectively.
In addition, when \(i = 2\), we read the tuples \(\boldsymbol{Y}\) and \(\boldsymbol{y}\) of Lemma~\ref{lem:relative-regulator-polynomial} as \(\boldsymbol{Z}\) and \(\boldsymbol{z}\), respectively.
Then there exist nonzero multihomogeneous polynomials
\[
P_{1} \in \QQ[\boldsymbol{X},\boldsymbol{Y}],
\qquad P_{2} \in \QQ[\boldsymbol{X},\boldsymbol{Z}]
\]
with respect to \((\boldsymbol{X},\boldsymbol{Y})\) and \((\boldsymbol{X},\boldsymbol{Z})\), respectively, such that
\[
\reg(K_{1}) = |P_{1}(\boldsymbol{x},\boldsymbol{y})|,
\qquad \reg(K_{2}) = |P_{2}(\boldsymbol{x},\boldsymbol{z})|.
\]
Since multiplication by \(-1\) preserves both the property of being nonzero and the multidegree of a polynomial, we may assume, after changing the signs of \(P_{1}\) and \(P_{2}\) if necessary, that
\[
\reg(K_{1}) = P_{1}(\boldsymbol{x},\boldsymbol{y}),
\qquad \reg(K_{2}) = P_{2}(\boldsymbol{x},\boldsymbol{z}).
\]
Furthermore, if we regard \(P_{1}\) and \(P_{2}\) as elements of \(\QQ[\boldsymbol{X},\boldsymbol{Y},\boldsymbol{Z}]\),
then the multidegree of \(P_{1}\) with respect to \((\boldsymbol{X},\boldsymbol{Y},\boldsymbol{Z})\) is
\[
(\delta, \rho_{1}, 0),
\]
and the multidegree of \(P_{2}\) is
\[
(\delta, 0, \rho_{2}).
\]
In particular, \(P_{1}\) is nonzero and each of its monomials has total degree \(\delta + \rho_{1} = r(K_{1}) > 0\);
hence at least one indeterminate appears in \(P_{1}\).
Thus \(\delta_{+} + e_{1} \ge 1\), and therefore \(\delta_{+} + e_{1} + e_{2} \ge 1\).
It follows from Lemma~\ref{lem:relative-log-independence} that
all components of \(\boldsymbol{x}, \boldsymbol{y}, \boldsymbol{z}\) are linearly independent over \(\QQ\).
Since these components are logarithms of algebraic numbers that are positive real numbers, WSC implies that they are algebraically independent over \(\QQ\).

By inequalities~\eqref{eq:multidegree-inequalities} and Lemma~\ref{lem:shared-homogeneous-pair}, \(P_{1}\) and \(P_{2}\) are algebraically independent over \(\QQ\).
Hence \(H(P_{1},P_{2})\) is a nonzero polynomial for every nonzero polynomial
\[
H \in \QQ[T_{1},T_{2}].
\]
Since all components of \(\boldsymbol{x},\boldsymbol{y},\boldsymbol{z}\) are algebraically independent over \(\QQ\), we obtain
\[
H\bigl(\reg(K_{1}),\reg(K_{2})\bigr) = H\bigl(P_{1}(\boldsymbol{x},\boldsymbol{y}),P_{2}(\boldsymbol{x},\boldsymbol{z})\bigr) \ne 0.
\]
Since \(H\) was arbitrary, \(\reg(K_{1})\) and \(\reg(K_{2})\) are algebraically independent over \(\QQ\).

Next, we prove that \textup{(ii)} implies \textup{(i)}.
Assume \textup{(ii)}.

If \(r(K_{i_{0}}) = 0\) for some \(i_{0} \in \{1, 2\}\), then \(\reg(K_{i_{0}}) = 1\) by definition,
and the nonzero polynomial \(T_{i_{0}} - 1 \in \QQ[T_{1}, T_{2}]\) vanishes at the point \(\bigl(\reg(K_{1}), \reg(K_{2})\bigr)\), which contradicts \textup{(ii)}.
Hence \(r(K_{1}) > 0\) and \(r(K_{2}) > 0\).

Moreover, if \(K_{1}^{+} = K_{2}^{+}\), then Lemma~\ref{lem:same-real-part-ratio} gives
\[
q \coloneqq \frac{\reg(K_{1})}{\reg(K_{2})} \in \QQ_{ > 0},
\]
and the point \(\bigl(\reg(K_{1}), \reg(K_{2})\bigr)\) is a zero of the nonzero polynomial \(T_{1} - qT_{2} \in \QQ[T_{1}, T_{2}]\), which again contradicts \textup{(ii)}.
Thus \(K_{1}^{+} \ne K_{2}^{+}\).

Therefore, \textup{(ii)} implies \textup{(i)}.
\end{proof}

\begin{proof}[Proof of Corollary~\ref{cor:ratio-transcendence}]
We first prove that \textup{(i)} implies \textup{(ii)}.
Assume that \(K_{1}^{+} \ne K_{2}^{+}\); we show that \(\reg(K_{1})/\reg(K_{2})\) is a transcendental number.

The assumption \(K_{1}^{+} \ne K_{2}^{+}\) is invariant under the interchange of the indices \(1\) and \(2\),
and a positive real number is a transcendental number if and only if its reciprocal is a transcendental number.
Hence, after interchanging \(K_{1}\) and \(K_{2}\) if necessary, we may assume that \(r(K_{1}) \ge r(K_{2})\).

Suppose first that \(r(K_{2}) > 0\). Then \(r(K_{1}) \ge r(K_{2}) > 0\),
and Theorem~\ref{thm:galois-pair-criterion} shows that \(\reg(K_{1})\) and \(\reg(K_{2})\) are algebraically independent over \(\QQ\).
Since \(\QQalg/\QQ\) is an algebraic extension, these regulators are also algebraically independent over \(\QQalg\)~\cite[Theorem~4.3.5]{Roman2006}.
Put \(\alpha \coloneqq \reg(K_{1})/\reg(K_{2})\), and assume that \(\alpha\) is an algebraic number.
Let \(T_{1}, T_{2}\) be indeterminates. Since \(\alpha \in \QQalg\), we have
\(T_{1} - \alpha T_{2} \in \QQalg[T_{1}, T_{2}]\), and this polynomial is nonzero because the coefficient of \(T_{1}\) is \(1\).
On the other hand, the definition of \(\alpha\) gives
\[
\reg(K_{1}) - \alpha\reg(K_{2}) = 0,
\]
which contradicts the algebraic independence of \(\reg(K_{1})\) and \(\reg(K_{2})\) over \(\QQalg\).
Hence, in this case, \(\reg(K_{1})/\reg(K_{2})\) is a transcendental number.

Suppose next that \(r(K_{2}) = 0\). Then Lemma~\ref{lem:rank-zero}~\textup{(ii)} gives \(K_{2}^{+} = \QQ\).
The assumption \(K_{1}^{+} \ne K_{2}^{+}\) gives \(K_{1}^{+} \ne \QQ\), and thus Lemma~\ref{lem:rank-zero}~\textup{(ii)} again gives \(r(K_{1}) > 0\).
By the case \(m = 1\) of Theorem~\ref{thm:main-algebraic-independence}, \(\reg(K_{1})\) is a transcendental number.
Moreover, \(r(K_{2}) = 0\) gives \(\reg(K_{2}) = 1\), and therefore
\[
\frac{\reg(K_{1})}{\reg(K_{2})} = \reg(K_{1})
\]
is a transcendental number.

In either case, \(\reg(K_{1})/\reg(K_{2})\) is a transcendental number, and consequently \textup{(i)} implies \textup{(ii)}.

Next, \textup{(ii)} clearly implies \textup{(iii)} because every rational number is an algebraic number.

Finally, we prove that \textup{(iii)} implies \textup{(i)}.
Assume \textup{(iii)}.
If \(K_{1}^{+} = K_{2}^{+}\), then Lemma~\ref{lem:same-real-part-ratio} gives
\[
\frac{\reg(K_{1})}{\reg(K_{2})} \in \QQ_{ > 0} \subseteq \QQ,
\]
which contradicts \textup{(iii)}.
Hence \(K_{1}^{+} \ne K_{2}^{+}\).
Therefore, \textup{(iii)} implies \textup{(i)}.

We now prove the last assertion of the corollary.
By the equivalence of \textup{(ii)} and \textup{(iii)}, \(\reg(K_{1})/\reg(K_{2})\) is a transcendental number unless it is a rational number.
In addition, a rational number is not a transcendental number.
Thus the ratio \(\reg(K_{1})/\reg(K_{2})\) is either a rational number or a transcendental number.
The statement that this ratio is a rational number if and only if \(K_{1}^{+} = K_{2}^{+}\) is merely a restatement of the equivalence of \textup{(i)} and \textup{(iii)}.
\end{proof}

\begin{proof}[Proof of Corollary~\ref{cor:totally-real-pair}]
For each \(i \in \{1, 2\}\), \(K_{i}\) is a totally real number field different from \(\QQ\), and hence
\[
r(K_{i}) = [K_{i}:\QQ] - 1 > 0.
\]
Moreover, since \(K_{i}\) is totally real,
\[
K_{i}^{+} = K_{i}.
\]
Therefore, condition \textup{(i)} of Theorem~\ref{thm:galois-pair-criterion} is equivalent to \(K_{1} \ne K_{2}\).
Hence Theorem~\ref{thm:galois-pair-criterion} yields the conclusion.
\end{proof}

\section{On the assumptions of the main theorems}\label{sec:remarks}

In this section, we examine the assumptions of the two main theorems and of their corollaries through concrete examples.

\subsection{The Galois assumption}\label{subsec:galois-hypothesis}

The following example shows that, if we remove the Galois assumption from the two main theorems and their corollaries, then the respective conclusions do not hold in general.

\begin{example}\label{ex:non-galois}
Following Cruz~\cite[Example~6.2]{Cruz2025}, we use two septic fields that originate from Mantilla-Soler~\cite[Remark~2.11]{MantillaSoler2015}.
That is, let \(T\) be an indeterminate, and define \(f_{1}, f_{2} \in \QQ[T]\) by
\[
\begin{aligned}
f_{1}(T)
& \coloneqq T^{7} - 2T^{6} - 47T^{5} + 25T^{4} + 755T^{3} + 496T^{2} - 3782T - 5217,\\
f_{2}(T)
& \coloneqq T^{7} - 2T^{6} - 47T^{5} - 8T^{4} + 480T^{3} + 793T^{2} + 233T + 19.
\end{aligned}
\]
The polynomials \(f_{1}, f_{2}\) are irreducible over \(\QQ\), and all roots of each of them are real.
For each \(i \in \{1, 2\}\), fix a real root \(\theta_{i} \in \RR\) of \(f_{i}\) and put \(K_{i} \coloneqq \QQ(\theta_{i})\).
The isomorphism class, the rank of the unit group, and the regulator of \(K_{i}\) do not depend on the choice of the root \(\theta_{i}\).
These fields are non-isomorphic totally real septic fields, and they are arithmetically equivalent~\cite[Remark~2.11]{MantillaSoler2015}.
Furthermore, \cite[Example~6.2]{Cruz2025} shows that
\[
\reg(K_{1}) = \reg(K_{2}).
\]
As stated in Remark~\ref{rem:regulator-gram-determinant}, the regulator of that paper gives the same value as \(\reg\) of the present paper.
As stated in Section~\ref{sec:introduction}, the product of the class number and the regulator is the same for two arithmetically equivalent number fields;
hence the above equality is equivalent to the coincidence of the class numbers of \(K_{1}\) and \(K_{2}\).

Since both \(K_{1}\) and \(K_{2}\) are totally real septic fields, we obtain
\[
r(K_{1}) = r(K_{2}) = 6 > 0.
\]
In addition, by the multiplicativity of degrees, \([K_{1} \cap K_{2}:\QQ]\) is a divisor of \(7\).
Since \(7\) is a prime number, \([K_{1} \cap K_{2}:\QQ]\) is \(1\) or \(7\).
If \([K_{1} \cap K_{2}:\QQ] = 7\), then
\[
[K_{i}:K_{1} \cap K_{2}] = 1
\]
for each \(i \in \{1, 2\}\), and hence \(K_{1} = K_{1} \cap K_{2} = K_{2}\), which contradicts the fact that the two fields are non-isomorphic.
Therefore, \([K_{1} \cap K_{2}:\QQ] = 1\), and
\[
K_{1} \cap K_{2} = \QQ.
\]

Furthermore, if one of two arithmetically equivalent number fields is Galois over \(\QQ\), then the two fields are isomorphic~\cite{Perlis1977}.
Since \(K_{1}\) and \(K_{2}\) are arithmetically equivalent but non-isomorphic, neither of them is Galois over \(\QQ\).

The unit groups of \(K_{1}\) and \(K_{2}\) both have positive rank, and the unit group of \(K_{1} \cap K_{2}\) has rank \(0\).
Therefore, if we remove only the Galois assumption from the assumptions of Theorem~\ref{thm:main-algebraic-independence},
then the conclusion of algebraic independence in that theorem does not hold in general.
Note that, in this example, the two fields are not even isomorphic although their regulators coincide.

Moreover, since \(K_{1}\) and \(K_{2}\) are totally real, we have \(K_{i}^{+} = K_{i}\) \((i = 1, 2)\), and since the two fields are non-isomorphic, we have \(K_{1}^{+} \ne K_{2}^{+}\).
On the other hand, \(\reg(K_{1})/\reg(K_{2}) = 1\).
Thus, in both Theorem~\ref{thm:galois-pair-criterion} and Corollary~\ref{cor:ratio-transcendence}, condition \textup{(i)} holds but condition \textup{(ii)} does not.
Consequently, we cannot remove the Galois assumption from Theorem~\ref{thm:galois-pair-criterion}, Corollary~\ref{cor:ratio-transcendence}, or Corollary~\ref{cor:totally-real-pair} either.
\end{example}

\subsection{The assumption on the rank of the unit group}\label{subsec:positive-unit-rank}

Take an arbitrary integer \(m \ge 1\) and arbitrary number fields \(K_{1}, \dots, K_{m}\).
If \(r(K_{i_{0}}) = 0\) for some \(i_{0} \in \{1, \dots, m\}\), then \(\reg(K_{i_{0}}) = 1\) by definition.
In this case, if \(T_{1}, \dots, T_{m}\) are indeterminates, then
the nonzero polynomial \(T_{i_{0}} - 1 \in \QQ[T_{1}, \dots, T_{m}]\) vanishes at the point \(\bigl(\reg(K_{1}), \dots, \reg(K_{m})\bigr)\).
Therefore, \(\reg(K_{1}), \dots, \reg(K_{m})\) are not algebraically independent over \(\QQ\).
In other words, if these regulators are algebraically independent over \(\QQ\), then the unit groups of \(K_{1}, \dots, K_{m}\) all have positive rank.
This argument uses neither the Galois property nor WSC.

In contrast, Corollary~\ref{cor:ratio-transcendence} also covers the case where the unit group of one or both of the fields has rank \(0\).

\begin{example}\label{ex:rank-zero-ratio}
Consider the two imaginary quadratic fields
\[
\begin{aligned}
K_{1}
& \coloneqq \QQ(\sqrt{-1}),\\
K_{2}
& \coloneqq \QQ(\sqrt{-2}).
\end{aligned}
\]
Both fields are Galois over \(\QQ\), and the maximal real subfield of each of them is \(\QQ\).
Moreover, since
\[
\reg(K_{1}) = \reg(K_{2}) = 1,
\]
we have
\[
\frac{\reg(K_{1})}{\reg(K_{2})} = 1.
\]
This example falls under the case \(K_{1}^{+} = K_{2}^{+}\), and the ratio of the regulators is a rational number, in accordance with the last assertion of Corollary~\ref{cor:ratio-transcendence}.
\end{example}

\subsection{The condition on the maximal real subfields of two fields}\label{subsec:maximal-real-subfields}

We cannot replace the condition \(K_{1}^{+} \ne K_{2}^{+}\) in Theorem~\ref{thm:galois-pair-criterion}
simply by the condition \(K_{1} \ne K_{2}\).

\begin{example}\label{ex:cm-intersection}
Consider the two fields
\[
\begin{aligned}
K_{1}
& \coloneqq \QQ(\sqrt{10}, \sqrt{-21}),\\
K_{2}
& \coloneqq \QQ(\sqrt{10}, \sqrt{-42}).
\end{aligned}
\]
These are distinct quartic CM fields that are Galois over \(\QQ\)~\cite[Example~7.1]{IizukaKonomi2026},
and \(r(K_{i}) = 1\) \((i = 1, 2)\).

Both fields contain the real quadratic field \(\QQ(\sqrt{10})\), and each of them is a quadratic extension of \(\QQ(\sqrt{10})\).
Hence
\[
\QQ(\sqrt{10}) \subseteq K_{1} \cap K_{2} \subseteq K_{1}.
\]
Moreover, since \([K_{1}:\QQ(\sqrt{10})] = 2\), the intersection \(K_{1} \cap K_{2}\) is either \(\QQ(\sqrt{10})\) or \(K_{1}\).
In the latter case, \(K_{1} \subseteq K_{2}\), and \([K_{1}:\QQ] = [K_{2}:\QQ]\) gives \(K_{1} = K_{2}\), which contradicts the fact that the two fields are distinct. Therefore,
\[
K_{1} \cap K_{2} = \QQ(\sqrt{10}) \ne \QQ.
\]
In particular,
\[
r(K_{1} \cap K_{2}) = 1 > 0.
\]
In addition, by the result of Iizuka and Konomi~\cite[Example~7.1]{IizukaKonomi2026}, we have
\[
\reg(K_{1}) = \reg(K_{2}) = 2\log(3 + \sqrt{10}).
\]

Furthermore, since both fields are quadratic extensions of the real quadratic field \(\QQ(\sqrt{10})\) and are not contained in \(\RR\),
the maximal real subfields of the two fields are
\[
K_{1}^{+} = K_{2}^{+} = \QQ(\sqrt{10}).
\]

This example shows that we cannot replace the condition \(K_{1}^{+} \ne K_{2}^{+}\) in Theorem~\ref{thm:galois-pair-criterion} by \(K_{1} \ne K_{2}\).
Moreover, this example falls under the case \(K_{1}^{+} = K_{2}^{+}\), and the ratio of the regulators
\[
\frac{\reg(K_{1})}{\reg(K_{2})} = 1
\]
is indeed a rational number, in accordance with the last assertion of Corollary~\ref{cor:ratio-transcendence}.

In contrast, if two distinct totally real number fields that are Galois over \(\QQ\) are both different from \(\QQ\),
then, by Corollary~\ref{cor:totally-real-pair}, their regulators are algebraically independent over \(\QQ\) under WSC.
\end{example}

\subsection{The assumption on the intersections in a finite family}\label{subsec:intersection-hypothesis}

The following example shows that we cannot replace the assumption on the intersections in Theorem~\ref{thm:main-algebraic-independence}
by the condition that the maximal real subfields are pairwise distinct.

\begin{example}\label{ex:distinct-maximal-real-subfields}
Consider the four totally real number fields
\[
\begin{aligned}
K_{1}
& \coloneqq \QQ(\sqrt{2}),\\
K_{2}
& \coloneqq \QQ(\sqrt{3}),\\
K_{3}
& \coloneqq \QQ(\sqrt{6}),\\
K_{4}
& \coloneqq \QQ(\sqrt{2}, \sqrt{3}).
\end{aligned}
\]
All of these fields are Galois over \(\QQ\), and they are pairwise distinct.
Moreover, since they are all totally real, we have
\[
K_{i}^{+} = K_{i} \quad (1 \le i \le 4),
\]
and the maximal real subfields \(K_{1}^{+}, K_{2}^{+}, K_{3}^{+}, K_{4}^{+}\) are pairwise distinct.
Furthermore,
\[
r(K_{1}) = r(K_{2}) = r(K_{3}) = 1,
\qquad r(K_{4}) = 3.
\]

For each \(i \in \{1,2,3\}\), take a fundamental unit \(\varepsilon_{i}\) of \(K_{i}\) such that \(\varepsilon_{i} > 1\).
The definition of the regulator gives
\[
\reg(K_{i}) = \log\varepsilon_{i} \quad (1 \le i \le 3).
\]

For each \(a,b \in \{\pm 1\}\), define the real embedding \(\sigma_{a,b}\) of \(K_{4}\) by
\[
\sigma_{a,b}(\sqrt{2}) = a\sqrt{2},
\qquad \sigma_{a,b}(\sqrt{3}) = b\sqrt{3}.
\]
Since \(\Norm_{K_{i}/\QQ}(\varepsilon_{i}) \in \{\pm 1\}\) for each \(i \in \{1,2,3\}\), we obtain
\[
\begin{aligned}
\log|\sigma_{a,b}(\varepsilon_{1})|
& = a\reg(K_{1}),\\
\log|\sigma_{a,b}(\varepsilon_{2})|
& = b\reg(K_{2}),\\
\log|\sigma_{a,b}(\varepsilon_{3})|
& = ab\reg(K_{3})
\end{aligned}
\]
for all \(a,b \in \{\pm 1\}\).

We regard \(\varepsilon_{1}, \varepsilon_{2}, \varepsilon_{3}\) as units of \(K_{4}\),
arrange the real embeddings of \(K_{4}\) in the order \(\sigma_{1,1}, \sigma_{-1,1}, \sigma_{1,-1}, \sigma_{-1,-1}\), and define the logarithmic matrix \(\mathcal{M}_{K_{4}}(\varepsilon_{1}, \varepsilon_{2}, \varepsilon_{3})\).
The last row of this matrix corresponds to \(\sigma_{-1,-1}\).
By the above equalities, the matrix \(\mathcal{M}_{K_{4}}^{(4)}(\varepsilon_{1}, \varepsilon_{2}, \varepsilon_{3})\) obtained by removing this row is
\[
\begin{pmatrix}
\reg(K_{1})
& \reg(K_{2})
& \reg(K_{3})\\
-\reg(K_{1})
& \reg(K_{2})
& -\reg(K_{3})\\
\reg(K_{1})
& -\reg(K_{2})
& -\reg(K_{3})
\end{pmatrix}.
\]
The absolute value of its determinant is
\[
4\reg(K_{1})\reg(K_{2})\reg(K_{3}) > 0.
\]
Now take arbitrary \(A_{1}, A_{2}, A_{3} \in \ZZ\), and assume that
\[
\varepsilon_{1}^{A_{1}}\varepsilon_{2}^{A_{2}}\varepsilon_{3}^{A_{3}} \in \mu(K_{4}).
\]
Every root of unity has absolute value \(1\) under every embedding.
Accordingly, if we apply \(\sigma_{1,1}, \sigma_{-1,1}, \sigma_{1,-1}\), respectively, to the above product
and take the logarithms of the absolute values of the resulting values, then we obtain
\[
\begin{pmatrix}
\reg(K_{1})
& \reg(K_{2})
& \reg(K_{3})\\
-\reg(K_{1})
& \reg(K_{2})
& -\reg(K_{3})\\
\reg(K_{1})
& -\reg(K_{2})
& -\reg(K_{3})
\end{pmatrix}
\begin{pmatrix}
A_{1}\\
A_{2}\\
A_{3}
\end{pmatrix}
=
\begin{pmatrix}
0\\
0\\
0
\end{pmatrix}.
\]
Since the determinant of the \(3 \times 3\) matrix on the left-hand side is nonzero, we obtain
\[
A_{1} = A_{2} = A_{3} = 0.
\]
Therefore, the cosets
\[
\varepsilon_{1}\mu(K_{4}),
\qquad \varepsilon_{2}\mu(K_{4}),
\qquad \varepsilon_{3}\mu(K_{4})
\]
are multiplicatively independent in \(\mathcal{O}_{K_{4}}^{\times}/\mu(K_{4})\).
In addition, since
\[
r(K_{4}) = 3,
\]
the subgroup generated by these cosets has finite index in \(\mathcal{O}_{K_{4}}^{\times}/\mu(K_{4})\).
Let
\[
I \coloneqq \left[ \mathcal{O}_{K_{4}}^{\times}/\mu(K_{4}) : \left\langle \varepsilon_{1}\mu(K_{4}), \varepsilon_{2}\mu(K_{4}), \varepsilon_{3}\mu(K_{4}) \right\rangle \right] \in \ZZ_{ > 0}
\]
be its index.

The index formula~\eqref{eq:index-formula} for logarithmic matrices gives
\[
I\reg(K_{4}) = 4\reg(K_{1})\reg(K_{2})\reg(K_{3}).
\]
Hence, if \(T_{1}, T_{2}, T_{3}, T_{4}\) are indeterminates, then the nonzero polynomial
\[
IT_{4} - 4T_{1}T_{2}T_{3} \in \QQ[T_{1}, T_{2}, T_{3}, T_{4}]
\]
vanishes at the point \(\bigl(\reg(K_{1}), \reg(K_{2}), \reg(K_{3}), \reg(K_{4})\bigr)\).
Thus the regulators
\[
\reg(K_{1}), \reg(K_{2}), \reg(K_{3}), \reg(K_{4})
\]
are not algebraically independent over \(\QQ\).

Consequently, we cannot simply replace the assumption on the intersections in Theorem~\ref{thm:main-algebraic-independence}
by the condition that the maximal real subfields are pairwise distinct.

On the other hand, the above four fields are all totally real number fields different from \(\QQ\), are Galois over \(\QQ\), and are pairwise distinct.
Hence, by Corollary~\ref{cor:totally-real-pair}, any two of these regulators are algebraically independent over \(\QQ\) under WSC.
However, the four regulators together satisfy the above polynomial relation and are therefore not algebraically independent over \(\QQ\).
This polynomial relation holds without assuming WSC.

Moreover, we have
\[
K_{1} \cap K_{4} = K_{1},
\qquad r(K_{1} \cap K_{4}) = 1 > 0,
\]
but, as stated above, \(\reg(K_{1})\) and \(\reg(K_{4})\) are algebraically independent over \(\QQ\) under WSC.
Therefore, under WSC, the condition that the unit group of the intersection of two fields has rank \(0\) is not a necessary condition for the regulators of the two fields to be algebraically independent over \(\QQ\).

Furthermore, as stated above, \(K_{1}, K_{2}, K_{3}\) are Galois over \(\QQ\), and their unit groups all have positive rank.
In addition, as stated in Remark~\ref{rem:linear-independence-small-m}, we have \(K_{i} \cap K_{j} = \QQ\) for all distinct \(i, j \in \{1, 2, 3\}\), and hence \(r(K_{i} \cap K_{j}) = r(\QQ) = 0\).
Thus \(K_{1}, K_{2}, K_{3}\) satisfy the assumptions of Theorem~\ref{thm:main-algebraic-independence}, and \(\reg(K_{1}), \reg(K_{2}), \reg(K_{3})\) are algebraically independent over \(\QQ\) under WSC.
On the other hand, since \(K_{1} \cap (K_{2}K_{3}) = K_{1}\) as stated in Remark~\ref{rem:linear-independence-small-m}, we find
\[
r\bigl(K_{1} \cap (K_{2}K_{3})\bigr) = r(K_{1}) = 1 > 0.
\]
In this way, the assumption on the intersections in Theorem~\ref{thm:main-algebraic-independence} is a condition only on the intersections of two distinct fields,
and it does not require that the unit group of the intersection of one field with the compositum of the remaining fields have rank \(0\).
\end{example}

\section*{Acknowledgements}

The authors thank Professor Shin Nakano for helpful discussions.

\section*{Declaration of generative AI and AI-assisted technologies in the manuscript preparation process}

During the preparation of this paper, the authors used two generative AI tools, OpenAI's GPT-6 Astra and Anthropic's Claude Fable 5.2, in order to translate the Japanese manuscript into English, to proofread the text, to check the consistency of notation and terminology, to detect errors and gaps in the arguments, to look for examples and counterexamples, and to search the literature for prior and related work.
After using these tools, the authors carefully reviewed the output and decided what to adopt.
In addition, they themselves verified that every step of the arguments is correct and that the hypotheses of the theorems invoked are satisfied.
They also checked that the cited statements agree with the original sources and that the bibliographic information is accurate.
Every proof in this paper is complete on the basis of the arguments given in the text and the cited literature, and no step relies on the output of a generative AI tool.
The generative AI tools are not authors of this paper, and the authors take full responsibility for its content.


\begin{thebibliography}{99}

\bibitem{BiluBugeaudMignotte2014}
Y.~F.~Bilu, Y.~Bugeaud and M.~Mignotte,
\emph{The Problem of Catalan},
Springer, Cham, 2014.

\bibitem{Cruz2025}
J.~Cruz,
Classification of totally real number fields via their zeta function, regulator, and log-unit lattice,
preprint, arXiv:2512.03432v2, 2025.

\bibitem{GunMurtyRath2011}
S.~Gun, M.~Ram Murty and P.~Rath,
Transcendental nature of special values of \(L\)-functions,
\emph{Canad.~J.~Math.} \textbf{63} (2011), no.~1, 136--152.

\bibitem{IizukaKonomi2026}
Y.~Iizuka and Y.~Konomi,
An infinite family of pairs of distinct quartic Galois CM-fields
with the same discriminant and regulator,
\emph{Tokyo J.~Math.} \textbf{49} (2026), no.~1, 225--237.

\bibitem{Keune2023}
F.~Keune,
\emph{Number Fields},
Radboud University Press, Nijmegen, 2023.

\bibitem{Khanduja2022}
S.~K.~Khanduja,
\emph{A Textbook of Algebraic Number Theory},
UNITEXT, vol.~135,
Springer, Singapore, 2022.

\bibitem{Lang2002}
S.~Lang,
\emph{Algebra}, revised third edition,
Graduate Texts in Mathematics, vol.~211,
Springer-Verlag, New York, 2002.

\bibitem{MantillaSoler2015}
G.~Mantilla-Soler,
Weak arithmetic equivalence,
\emph{Canad.~Math.~Bull.} \textbf{58} (2015), no.~1, 115--127.

\bibitem{MilneANT}
J.~S.~Milne,
\emph{Algebraic Number Theory}, Version 3.08, 2020,
available at \url{https://www.jmilne.org/math/}.

\bibitem{MilneFT}
J.~S.~Milne,
\emph{Fields and Galois Theory}, Version 5.10, 2022,
available at \url{https://www.jmilne.org/math/}.

\bibitem{Perlis1977}
R.~Perlis,
On the equation \(\zeta_{K}(s) = \zeta_{K'}(s)\),
\emph{J.~Number Theory} \textbf{9} (1977), no.~3, 342--360.

\bibitem{Prasad2017}
D.~Prasad,
A refined notion of arithmetically equivalent number fields,
and curves with isomorphic Jacobians,
\emph{Adv.~Math.} \textbf{312} (2017), 198--208.

\bibitem{Roman2006}
S.~Roman,
\emph{Field Theory}, second edition,
Graduate Texts in Mathematics, vol.~158,
Springer, New York, 2006.

\bibitem{Serre2010}
D.~Serre,
\emph{Matrices: Theory and Applications}, second edition,
Graduate Texts in Mathematics, vol.~216,
Springer, New York, 2010.

\bibitem{SmitPerlis1994}
B.~de Smit and R.~Perlis,
Zeta functions do not determine class numbers,
\emph{Bull.~Amer.~Math.~Soc. (N.S.)} \textbf{31} (1994), no.~2, 213--215.

\bibitem{Sutherland2021}
A.~V.~Sutherland,
Stronger arithmetic equivalence,
\emph{Discrete Anal.} (2021), Paper No.~23, 23 pp.

\bibitem{Washington1997}
L.~C.~Washington,
\emph{Introduction to Cyclotomic Fields}, second edition,
Graduate Texts in Mathematics, vol.~83,
Springer-Verlag, New York, 1997.

\end{thebibliography}
\end{document}